\documentclass [11 pt , a4 paper , twoside, openright, BCOR5mm, cleardoublepage = empty]{article}

\input{Pack.sty}
\usepackage{color}
\usepackage{authblk}
\usepackage{subfigure}
 \usepackage{multirow}
\begin{document}
\title{
Reduced Basis POD-Galerkin Method for Parametrized Optimal Control Problems in Environmental Marine Sciences and Engineering}

\author[$\sharp$]{Maria Strazzullo}
\author[$\sharp$]{Francesco Ballarin}
\author[$\flat$]{Renzo Mosetti}
\author[$\sharp$]{Gianluigi Rozza}
\affil[$\sharp$]{mathlab, Mathematics Area, SISSA, Via Bonomea 265, I-34136 Trieste, Italy}
\affil[$\flat$]{National Institute of Oceanography and Experimental Geophysics, Via Beirut 2, I-34151 Trieste, Italy}
\date{}                     %% if you don't need date to appear
\setcounter{Maxaffil}{0}
\renewcommand\Affilfont{\itshape\small}
\maketitle

\begin{abstract}
In this work we propose reduced order methods as a suitable approach to face parametrized optimal control problems governed by partial differential equations, with applications in environmental marine sciences and engineering. Environmental parametrized optimal control problems are usually studied for different configurations described by several physical and/or geometrical parameters representing different phenomena and structures. The solution of parametrized problems requires a demanding computational effort. In order to save computational time, we rely on reduced basis techniques as a suitable and rapid tool to solve parametrized problems. We introduce general parametrized linear quadratic optimal control problems, and the saddle-point structure of their optimality system. Then, we propose a POD-Galerkin reduction of the optimality system. We test the resulting method on two environmental applications: a pollutant control in the Gulf of Trieste, Italy and a solution tracking governed by quasi-geostrophic equations describing North Atlantic Ocean dynamic. The two experiments underline how reduced order methods are a reliable and convenient tool to manage several environmental optimal control problems, for different mathematical models, geographical scale as well as physical meaning. The quasi-geostrophic optimal control problem is also presented in its nonlinear version.
\end{abstract}

% REQUIRED
\vspace{.5cm}

% REQUIRED
\textbf{AMS}: 
49J20, 76N25, 35Q35

\section{Introduction}

Parametrized optimal control problems (OCP($\boldsymbol{\mu}$)s) governed by parametrized partial differential equations (PDE($\boldsymbol{\mu}$)s) are usually complex and demanding, computationally speaking. In this case the parameter $\boldsymbol{\mu} \in \Cal P \subset \mathbb R^d$ could represent physical or geometrical features. In order to study different configurations, a rapid and suitable approach based on reduced order models could allow to face OCP($\boldsymbol{\mu}$)s in a low dimensional framework. Computational methods for OCP($\boldsymbol{\mu}$)s are a quite widespread tool in many contexts and fields: in shape optimization (see e.g. \cite{ delfour2011shapes, makinen, mohammadi2010applied}), in flow control (see e.g. \cite{de2007optimal, dede2007optimal, negri2015reduced, optimal}), in environmental applications (see e.g. \cite{phd, quarteroni2005numerical, quarteroni2007reduced}). Reduced order methods are a strategy that decreases the computational costs and simplifies the resolution of simulations governed by partial differential equations (see e.g. \cite{hesthaven2015certified,RozzaHuynhManzoni2013, RozzaHuynhPatera2008, prud2002reliable}). Reduction techniques has also been exploited in order to manage OCP($\pmb \mu$) in a faster way: reduced order methods for optimal control problems are theoretically and experimentally treated in several works different for state equations, methodology and purposes (see e.g. \cite{bader2016certified, Dedè2012, gerner2012certified, karcher2014certified,
kunisch1999control, kunisch2008proper, negri2013reduced, negri2015reduced, quarteroni2007reduced}). 
\\In this work, we focus on reduced order modelling for optimal control problems with quadratic cost functional constrained to linear and nonlinear PDE($\boldsymbol{\mu}$)s dealing with applications in environmental marine sciences. Optimal control theory fits well in these fields, as it is at the basis of forecasting models and it could be used to make previsions on several natural phenomena (see e.g. \cite{ghil1991data,kalnay2003atmospheric,tziperman1989optimal}). In order to describe different configurations, simulations have to be run for several values of $\boldsymbol{\mu}$. For this reason reduced order methods are presented as an important resource to manage these problems. The main novelty of this work deals with the use of %\reviewerA{well-known} 
reduced order modelling techniques in geographically realistic experiments involved in marine sciences with environmental purposes. Two examples will be discussed:
\begin{enumerate}
\item a first example dealing with a pollutant control in the Gulf of Trieste, Italy. In this case we aimed at underlining how parametrized optimal control problems could be useful in order to study several potential configurations describing different phenomena in this specific geographic area. These experiments could improve the monitoring strategies (presented, for example, in \cite{mosetti2005innovative, shiganova2009native}) of the marine environment of the Gulf of Trieste;
\item an Oceanographic solution tracking governed by quasi-geostrophic equations (see e.g. \cite{cavallini2012quasi, kim2015b, san2015stabilized}), describing the North Atlantic Ocean dynamics. This problem could be linked to \textit{data assimilation} based forecasting models (see e.g. \cite{behringer1998improved,carton2000simple, carton2008reanalysis, ghil1991data, kalnay2003atmospheric, tziperman1989optimal}) which are used to simulate possible scenarios and to analyse and predict climatological phenomena. The solution tracking is presented both in the linear and in the nonlinear version.
 \vspace{.2cm}
\end{enumerate} 

\no One of the main purposes of this work is to present reduced order methods as a reliable and useful tool to manage environmental simulations, both for problems characterized by large scales as well as small ones. Indeed, the two experiments that we have analysed are different under many aspects: equations used, scale analysis, geographical regions considered, dynamics involved. Despite that, they have in common a physical parametrized setting describing several configurations and modelling different natural phenomena. Many computational resources are needed to solve an optimal control problem in environmental sciences, most of all when parameter-dependent simulations have to be run many times in order to study and analyse several results, representing very different physical and natural aspects. Reduced order methods allow to recast a computational demanding problem, the \textit{``truth" problem}, into a new fast and reliable low-dimensional formulation, the \textit{reduced problem}. The latter is formulated as a Galerkin projection into reduced spaces, generated by basis functions chosen through a proper orthogonal decomposition sampling  algorithm, as presented in \cite{ballarin2015supremizer, burkardt2006pod, Chapelle2013, hesthaven2015certified}.
\\This work will show how convenient the reduced formulation is, since it is capable to give \textit{real-time} results, while environmental simulations based on classical approximation methods may take a very long time. The computational time saved for the reduced simulations could be invested in the study of many scenarios in order to achieve a deeper knowledge of ecological and climatological phenomena, that are very hard to be analyzed and understood, since they are linked to several natural aspects, as well as anthropic features.
\\ 
\no The paper is outlined as follows. In section 2, first of all the saddle-point structure for linear quadratic parametrized optimal control problems and its Finite Element (FE) approximation are briefly discussed (see \cite{tesi, fernandez2003control,gunzburger2003perspectives}). Section 3 aims at introducing the reduced order approximation for OCP($\pmb \mu$)s (following \cite{hesthaven2015certified, ito1998reduced, karcher2014certified}) and POD sampling algorithm for OCP($\boldsymbol{\mu}$)s with a brief mention of aggregated reduced space strategy (used in \cite{ dede2010reduced, negri2015reduced,negri2013reduced}) and affine decomposition (see e.g. \cite{hesthaven2015certified}). In section 4, numerical results dealing with reduced order methods applied to environmental marine \nl{linear quadratic OCP($\pmb \mu$)s }are detailed. \nl{Finally, in section 5, the nonlinear version of a solution tracking governed by quasi-geostrophic equation is presented. } Conclusions follow in section 6.

\section{\color{black}{Linear Quadratic OCP($\boldsymbol{\mu}$)s: Problem Formulation and Finite Element Approximation}}
\label{conty}

This second section aims at describing linear quadratic OCP($\boldsymbol{\mu}$)s exploiting their saddle-point formulation. This strategy is illustrated in many works and in many applications (see \cite{negri2015reduced, negri2013reduced, rozza2012reduction, schoberl2007symmetric}). Saddle-point formulation is advantageous, since the results about its well-posedness are well known in literature (see \cite{boffi2013mixed, brezzi1974existence}). Then, we will briefly focus on the Finite Element (FE) ``truth" approximation of OCP($\boldsymbol{\mu}$)
%\end{itemize}

%\subsection{\nl{General Problem Formulation: Lagrangian Representation}}
%\nl{First of all, let us discuss the general problem formulation}

%This section aims at describing linear quadratic OCP($\boldsymbol{\mu}$)s exploiting their saddle-point formulation. This strategy is illustrated in many works and in many applications (see \cite{negri2015reduced, negri2013reduced, rozza2012reduction, schoberl2007symmetric}). Saddle-point formulation is advantageous, since the results about its well-posedness are well known in literature (see \cite{boffi2013mixed, brezzi1974existence}). Then, we will briefly focus on the Finite Element (FE) ``truth" approximation of OCP($\boldsymbol{\mu}$)s.
\subsection{Linear Quadratic OCP($\boldsymbol{\mu}$)s and Saddle-Point Structure}
\label{LQ}
In the treatment of linear quadratic OCP($\boldsymbol{\mu}$)s and their saddle-point structures we will essentially follow \cite{ tesi, negri2015reduced, negri2013reduced}. %\nl{Exploiting saddle-point formulation  is very advantageous, since the results about its well-posedness are well known in literature (see \cite{boffi2013mixed, brezzi1974existence});}
Let us consider $\Omega \subset \mathbb R^n$ an open and bounded domain with Lipschitz boundary 
$\partial \Omega$. Let $Y$ and $U$  be Hilbert spaces for state and control, respectively. Moreover, let $Z \supset Y$ be an Hilbert space, in which the observation is taken. Furthermore, let us define a compact set of parameters $\Cal P \subset \mathbb R^d$, for $d \geq 1$, where $\boldsymbol{\mu} = (\mu_1, \dots, \mu_d)$ is considered. 
We consider $Y$ as the \emph{adjoint} Hilbert space, in other words, in our applications, the adjoint and the state spaces will always coincide. Then, the linear constraint equation is defined by:
\begin{equation}
\label{eq : stato}
a(y,q; \boldsymbol{\mu}) = c(u,q; \boldsymbol{\mu}) + \la G(\boldsymbol{\mu}), q \ra \hspace{1cm} \forall q \in Y ,
\end{equation}
\vspace{-.05cm}
\no where $a : Y \times Y \rightarrow \mathbb R$ represents a continuous bilinear state operator, 
$c : U \times Y \rightarrow \mathbb R$ is a continuous bilinear form describing the role of the control in the formulation of the problem and $G \in\red{Y}\dual$ gathers information about forcing and boundary terms.
Set a constant $\alpha > 0$, the quadratic objective functional is given by:
\begin{equation}
\label{obj}
\displaystyle J(y,u; \boldsymbol{\mu}) = \half m(y -y_d(\boldsymbol{\mu}), y - y_d(\boldsymbol{\mu}); \boldsymbol{\mu}) + \alf n(u,u; \boldsymbol{\mu}),
\end{equation}
\no where $y_d(\boldsymbol{\mu}) \in Z$, $m : Z \times Z \rightarrow \mathbb R$ and $n : U \times U \rightarrow \mathbb R$ are continuous bilinear forms representing the objective on the state variable and a penalization for the control variable, respectively.
We remark that in writing the bilinear forms $a, c, m$ and $n$, the dependence on the parameter is sometimes understood.
\\  An OCP($\boldsymbol{\mu}$) problem can be formalized as follows: given $\boldsymbol{\mu}$, solve
\begin{equation}
\label{eq: controllo}
\underset{(y(\boldsymbol{\mu}),u(\boldsymbol{\mu})) \in Y \times U}{\text{min }}J(y(\boldsymbol{\mu}),u(\boldsymbol{\mu}); \boldsymbol{\mu}) \hspace {1cm} \text{such that } (y(\boldsymbol{\mu}),u(\boldsymbol{\mu})) \in Y \times U\text{ satisfies } \eqref{eq : stato}.
\end{equation}
\no In order to recast the problem \eqref{eq: controllo} in a saddle-point formulation, let us define $X = Y \times U$. Being $x = (y,u)$ and $w = (z,v)$ two elements of $X$, we can endow $X$ with the scalar product 
$(x, w)_X = (y,z)_Y + (u, v)_U$ and with the norm $\norm {\cdot}_X = \sqrt{\cd_X}$. Now let us consider
\begin{align*} & \Cal A : X \times X \rightarrow  \mathbb R &&
\Cal A (x, w; \boldsymbol{\mu}) = m(y,z; \boldsymbol{\mu}) + \alpha n(u,v; \boldsymbol{\mu}) \hspace{1cm} &&& \forall x , w \in X, \\
& \Cal B : X \times \red{Y} \rightarrow \mathbb R && \Cal B(w, q; \boldsymbol{\mu}) = a(z,q; \boldsymbol{\mu}) - c(v,q; \boldsymbol{\mu}) \hspace{1cm} &&& \forall w \in X, \; \forall q \in \red{Y},\\
& F(\boldsymbol{\mu}) \in X\dual && \la F(\boldsymbol{\mu}), w \ra = m(y_d(\boldsymbol{\mu}),z; \boldsymbol{\mu}) \hspace{1cm} &&& \forall w \in X.
\end{align*}
Thanks to these quantities, the following functional can be defined:
\begin{equation}
\label{lin_func}
\Cal J(x; \boldsymbol{\mu}) = \half \Cal A(x, x; \boldsymbol{\mu}) - \la F(\boldsymbol{\mu}), x \ra.
\end{equation}
\no In \cite{bochev2009least, negri2013reduced}, it is shown that minimizing $\Cal J(x; \boldsymbol{\mu})$ with respect to $x \in X$ is equivalent to minimize $J(y,u; \boldsymbol{\mu})$ with respect $(y,u)\in Y \times U$. So problem \eqref{eq: controllo} is equivalent to
%$$J(y,u; \boldsymbol{\mu}) = \Cal J(x) + \Cal M(y_d),$$
%where $\displaystyle \Cal M(y_d) = \half m(y_d, y_d)$  is a constant term that does not give any contribution to the minimization of $J\cd$. For these reasons, it is now possible give a new formulation to the problem \eqref{eq: controllo}: find
\begin{equation}
\label{eq: SPOPT}
\underset{x \in X}{\text{min }}\Cal J(x; \boldsymbol{\mu}) \hspace{1cm } \text{such that } \hspace{1cm} \Cal B (x, q; \boldsymbol{\mu}) =
\la G(\boldsymbol{\mu}) , q \ra \hspace{1cm} \forall q \in \red{Y}. 
\end{equation}
The constrained optimization problem \eqref{eq: SPOPT} can be recast into an unconstrained optimization problem by defining the Lagrangian functional $\Lg :X \times \red{Y} \rightarrow \mathbb R $ as
\begin{equation}
\label{daje}
\Lg (x, p; \boldsymbol{\mu}) = \Cal J(x; \boldsymbol{\mu}) + \Cal B(x, p; \boldsymbol{\mu}) - \la G(\boldsymbol{\mu}), p \ra,
\end{equation}
\no where $p \in \red{Y}$ is the \textit{adjoint variable}. Thanks to the continuity the forms $a, c, m$ and $n$, the operators $\Cal A, \Cal B$ and $F$ are bounded and then, the Lagrangian is G\^{a}teaux derivable and the minimization problem \eqref{eq: SPOPT} is equivalent to finding the saddle-point of the Lagrangian functional \eqref{daje} and this leads to the saddle-point formulation (see \cite{bochev2009least, schoberl2007symmetric} as references). 
\\The new formulation of the problem is: given $\boldsymbol{\mu} \in \Cal P$, find $(x(\boldsymbol{\mu}), p(\boldsymbol{\mu})) \in X \times \red{Y}$ such that
\begin{equation}
\label{ZETA}
\begin{cases}
\Cal A(x(\boldsymbol{\mu}), w; \boldsymbol{\mu}) + \Cal B(w,p(\boldsymbol{\mu}); \boldsymbol{\mu}) = \la F(\boldsymbol{\mu}), w \ra & \forall w \in X, \\
\Cal B (x(\boldsymbol{\mu}), q; \boldsymbol{\mu}) = \la G(\boldsymbol{\mu}), q \ra & \forall q \in \red{Y}.
\end{cases}
\end{equation}
\\ The existence and the uniqueness of the solution \red{is provided by considering the adjoint space coinciding with the state space and, then, } (see \cite{ bochev2009least, tesi, negri2013reduced}) by the fulfillment of the hypotheses of the Brezzi's theorem reported in \cite{boffi2013mixed}. 
%\\ \textcolor{red}{\textbf{inizio aggiunta}}
%\\In order to guarantee the well-posedness of the optimal control problem, the bilinear forms must verify:
%\begin{enumerate}[(i)]
%\item continuity of $\Cal A(\cdot, \cdot; \boldsymbol{\mu})$ over $X \times X$ for all $\boldsymbol{\mu} \in \Cal P$;
%\item coercivity of $\Cal A(\cdot, \cdot; \boldsymbol{\mu})$ over $K = \{w \in X \; : \: \Cal B(w, q; \boldsymbol{\mu}) = 0, \; \forall q \in Q\};$
%\item continuity of $\Cal B(\cdot, \cdot; \boldsymbol{\mu})$ over $X \times Q$ for all $\boldsymbol{\mu} \in \Cal P$;
%\item the fulfillment of the \textit{inf-sup condition} of $\Cal B(\cdot, \cdot; \boldsymbol{\mu})$ over $X \times Q$, i.e. there exists a constant $\beta_0>0$ such that
%\begin{equation}
%\label{infsupcont}
%\beta(\boldsymbol{\mu}) = \inf_{q \in Q} \sup_{w \in X} \frac{\Cal B(w,q; \boldsymbol{\mu})}{\norm{w}_X \norm{q}_Q} > \beta_0 \hspace{1cm} \forall \boldsymbol{\mu} \in \Cal P.
%\end{equation} 
%\end{enumerate}
%\no The conditions (i)-(iv), can be extended to the discrete ``truth" Finite Element approximation and to the \textit{reduced} approximation, that will be introduced at a later time among this work. 
%\\ \textcolor{red}{\textbf{fine aggiunta}}

\subsection{Finite Element ``Truth" Approximation of OCP($\boldsymbol{\mu}$)s}
\label{PROFE}
\no Let $\{ \Cal T \disc \}$ be a triangulation over $\Omega$. In a Finite Element approximation $Y \disc = Y \cap\Cal X^{\Cal N}_r$
and $U \disc = U  \cap \Cal X^{\Cal N}_r$, where
$$
 \Cal X^{\Cal N}_r = \{ v \disc \in C^0(\overline \Omega) \; : \; v \disc |_{K} \in \mbb P_r, \; \; \forall K \in \Cal T \disc \}. 
$$ 
\no and $\mbb P_r$ represents the space of polynomials of degree at most equal to $r$ and $K$ a triangle of $\Cal T^{\Cal N}$. To have Brezzi's hypotheses guaranteed, \red{we also assume that the discretized state and adjoint spaces are coinciding } (see \cite{ChenQuarteroniRozza2015, tesi, negri2013reduced}). Let us consider the discrete product space $X^{\Cal N} = Y \disc \times U \disc \subset X$. The Galerkin Finite Element discretization of the saddle-point problem \eqref{ZETA} is: given $\boldsymbol{\mu} \in \Cal P$, find $(x^{\Cal N}(\boldsymbol{\mu}), p^{\Cal N}(\boldsymbol{\mu})) \in X^{\Cal N} \times 
\red{Y}^ {\Cal N}$ such that 

\begin{equation}
\label{FEOCP}
\begin{cases}
 \Cal A(x^{\Cal N}(\boldsymbol{\mu}),w\disc; \boldsymbol{\mu}) + \Cal B (w\disc,p(\boldsymbol{\mu})\disc; \boldsymbol{\mu}) = \la F(\boldsymbol{\mu}), w\disc \ra & \forall v\disc \in X\disc, \\
\Cal B(x\disc(\boldsymbol{\mu}),q\disc; \boldsymbol{\mu})= \la G(\boldsymbol{\mu}), q\disc \ra & \forall q\disc \in \red{Y}\disc. 
\end{cases}
\end{equation}

\no Let us focus on the algebraic structure of the system associated to (\ref{FEOCP}). The dimension of $X\disc$ and $\red{Y} \disc$ are respectively indicated with $\Cal N_X$ and $\Cal N_{\red{Y}}$. Let us define the basis of the finite dimensional spaces $X \disc$ and $\red{Y} \disc$ respectively as:
$$
\{ \varphi_i\}_{i = 1}^{\Cal N_X} \hspace{0.2cm} \text{and} \hspace{0.2cm} \{ \psi_j\}_{j = 1}^{\Cal N_{\red{Y}}}.
$$

\no We now can rewrite the solution $(x \disc (\boldsymbol{\mu}), p \disc (\boldsymbol{\mu})) \in X\disc \times \red{Y}\disc$ as:
$$
\displaystyle \Bigg(x\disc (\boldsymbol{\mu})  = \sum_{i=1}^{\Cal N_X} x_i^{\boldsymbol{\mu}} \varphi _i, \;  p \disc (\boldsymbol{\mu}) = \sum_{j=1}^{\Cal N_{\red{Y}}} p_j^{\boldsymbol{\mu}} \psi _j \Bigg ).
$$
\no Let us define 
$ A(\boldsymbol{\mu}) \in \mbb R^{\Cal N_X \times \Cal N_X }, 
B(\boldsymbol{\mu}) \in \mbb R^{\Cal N_{\red{Y}} \times \Cal N_X }, \mbf F(\boldsymbol{\mu}) \in \mbb R^{\Cal N_X}$ and
$\mbf G(\boldsymbol{\mu}) \in \mbb R^{\Cal N_{\red{Y}}}$ as follows:
$$A_{ij}(\boldsymbol{\mu}) = \Cal A(\varphi_i, \varphi_j; \boldsymbol{\mu}), \hspace{0.3cm} B_{ml}(\boldsymbol{\mu}) = \Cal B( \varphi_l, \psi_m; \boldsymbol{\mu}), 
\hspace{0.3cm} \mbf F_k(\boldsymbol{\mu}) = \la F(\boldsymbol{\mu}), \varphi_k \ra, \hspace{0.3cm} \mbf G_s(\boldsymbol{\mu}) = \la G(\boldsymbol{\mu}), \psi_s \ra.$$

\no From those quantities, we can build the following linear system, with a block structure:
\begin{equation}
\label{sistlin}
\begin{pmatrix}
A(\boldsymbol{\mu}) & B^T(\boldsymbol{\mu}) \\
B(\boldsymbol{\mu}) & 0
\end{pmatrix}
\begin{pmatrix}
\mbf x^{\boldsymbol{\mu}} \\
\mbf p^{\boldsymbol{\mu}}
\end{pmatrix}
=
\begin{pmatrix}
\mbf F(\boldsymbol{\mu}) \\
\mbf G(\boldsymbol{\mu})
\end{pmatrix},
\end{equation}
\no where $(\mbf x^{\boldsymbol{\mu}})_i = x_i^{\boldsymbol{\mu}}$ and $(\mbf p^{\boldsymbol{\mu}})_j = p_j^{\boldsymbol{\mu}}$. The approach proposed is known as $\textit{optimize-then-discretize}$ (see \cite{fernandez2003control,gunzburger2003perspectives}). In all the applications presented, the linear system is solved through \textit{one-shot} method (see \cite{schulz2009one, taasan1991one}).

\section{Reduced Basis Methods for Parametrized Optimal Control Problems}
In this section reduced basis methods for OCP($\boldsymbol{\mu}$)s are described. First of all, the general idea of reduced order approximation for OCP($\boldsymbol{\mu}$)s is given as proposed in \cite{phd, negri2015reduced, negri2013reduced, karcher2014certified}. Then, proper orthogonal decomposition (POD, see \cite{ballarin2015supremizer, burkardt2006pod, Chapelle2013, hesthaven2015certified} as references) and the theory of aggregated spaces will be introduced with some considerations about the efficiency of the method through affinity assumption, as presented in \cite{hesthaven2015certified}.

\subsection{Problem Formulation and Solution Manifold}
In section \ref{LQ}, we have already affirmed that a linear quadratic OCP($\boldsymbol{\mu}$)s could be formulated as a saddle-point problem of the form \eqref{ZETA}. Let us recall that $x(\boldsymbol{\mu}) = (y(\boldsymbol{\mu}), u(\boldsymbol{\mu}))$. In several cases, the association $\boldsymbol{\mu} \rightarrow (x(\boldsymbol{\mu}), p(\boldsymbol{\mu})) \in X \times \red{Y}$ defines a smooth \textit{solution manifold} of the form:
$$
\Cal M = \{ (x(\boldsymbol{\mu}), p(\boldsymbol{\mu}))\;| \; \boldsymbol{\mu} \in \Cal P\}.
$$
\no When the \textit{full order} problem \eqref{FEOCP} is solved, one finds the \textit{approximated solution manifold}:
$$
\Cal M \disc = \{ (x\disc(\boldsymbol{\mu}), p\disc(\boldsymbol{\mu}))\;| \; \boldsymbol{\mu} \in \Cal P\}.
$$
\no Reduced basis methods aim at building a good approximation of $\Cal M \disc$ through linear combination of properly chosen \textit{snapshots} $x \disc (\boldsymbol{\mu})$ and $ p \disc(\boldsymbol{\mu})$, assuming that the \tit{approximated manifold} has a smooth dependence from $\boldsymbol{\mu}$. In other words, the reduced spaces are built with \textit{full order} solutions computed for suitable parameters in $\Cal P$. Let us suppose to have already built $X_N \subset X \disc \subset X$ and $\red{Y}_N \subset \red{Y} \disc \subset \red{Y}$ as reduced product space and reduced adjoint space, respectively (the reduced spaces will be specified in section \ref{aggaff}). Than, the \textit{reduced problem} is formulated as follows: given $\boldsymbol{\mu} \in \Cal P$, find $(x_N(\boldsymbol{\mu}), p_N(\boldsymbol{\mu})) \in X_N \times \red{Y}_N$ such that
\begin{equation}
\label{RBOCP}
\begin{cases}
\Cal A(x_N(\boldsymbol{\mu}),w_N; \boldsymbol{\mu}) + \Cal B(w_N, p_N(\boldsymbol{\mu}); \boldsymbol{\mu}) = \la F(\boldsymbol{\mu}), w_N \ra & \forall w_N \in X_N, \\
\Cal B(x_N(\boldsymbol{\mu}), q_N; \boldsymbol{\mu}) =  \la G(\boldsymbol{\mu}), q_N \ra & \forall q_N \in \red{Y}_N.
\end{cases}
\end{equation}

\subsection{POD Algorithm for OCP($\boldsymbol{\mu}$)s}
\label{PODsec}
Let us focus our attention on POD algorithm used as sampling procedure for the construction of the reduced bases, as treated in \cite{ballarin2015supremizer, burkardt2006pod, Chapelle2013, hesthaven2015certified}. The POD approach is more suitable than any greedy algorithm for the application proposed in sections \ref{OCEANO} and \ref{nl_sec}, where we cannot rely on a posteriori error estimator, since the state equation may not be coercive (see footnote \ref{coe}). In order to apply the POD, a discrete and finite dimensional subset $\Cal P_h \subset \Cal P$ is needed. For this specific set of parameters, the \tit{discrete solution manifold} is defined as:
$$
\Cal M \disc (\Cal P_h) = \{ ( x \disc(\boldsymbol{\mu}), p\disc(\boldsymbol{\mu})) \; | \; \boldsymbol{\mu} \in \Cal P_h\}.
$$

\no The cardinality of $\Cal M \disc (\Cal P_h)$ is $M = |\Cal P_h|$. Naturally it holds
$
\Cal M \disc (\Cal P_h) \subset \Cal M\disc
$
since $\Cal P_h \subset \Cal P$. If $\Cal P_h$ is fine enough, $\Cal M \disc (\Cal P_h) $ is a good approximation of the discrete manifold $\Cal M \disc$. From now on, we will refer to the set of the linear combinations of elements of $\Cal M \disc (\Cal P_h)$ as $\mathbb M$.
The algorithm of POD is based on two processes:
\begin{enumerate}
\item sampling the parameter space $\Cal P_h$ in order to compute the \emph{full order solutions} at selected parameters,
\item a compression phase, where one discards the redundant information, respectively for state, control and adjoint variables.
\end{enumerate}

\no The $N$-spaces resulting from the POD algorithm minimize the following quantities, respectively:
\begin{equation}
\label{crit}
\sqrt{\frac{1}{M}
\sum_{\boldsymbol{\mu} \in \Cal P_h} \underset{z_N \in Y_N}{\text{\red{min} }} \norm{y\disc(\boldsymbol{\mu}) - z_N}_Y^2},
\hspace{.2cm}
\sqrt{\frac{1}{M}
\sum_{\boldsymbol{\mu} \in \Cal P_h} \underset{v_N \in U_N}{\text{\red{min}}} \norm{u\disc(\boldsymbol{\mu}) - v_N}_U^2},
\hspace{.2cm}
\sqrt{\frac{1}{M}
\sum_{\boldsymbol{\mu} \in \Cal P_h} \underset{q_N \in \red{Y}_N}{\text{\red{min} }} \norm{p\disc(\boldsymbol{\mu}) - q_N}^2_{\red{Y}}}.
\end{equation} 
\no Let us introduce an ordering on the parameters $\boldsymbol{\mu}_1, \dots, \boldsymbol{\mu}_M \in \Cal P_h$. This induces an ordering on the \tit{full order} solutions $y\disc(\boldsymbol{\mu}_1), \dots, y\disc(\boldsymbol{\mu}_M)$, $u\disc(\boldsymbol{\mu}_1), \dots, u\disc(\boldsymbol{\mu}_M)$ and $p\disc(\boldsymbol{\mu}_1), \dots, p\disc(\boldsymbol{\mu}_M)$. In order to build the POD-spaces, we define the symmetric and linear operators:
\begin{align*}
& \mbf C_y: \mathbb M \rightarrow \mathbb M  \hspace{.2cm}&&
\mbf C_y(z\disc) = \frac{1}{M} \sum_{m=1}^{M} (z\disc, y\disc(\boldsymbol{\mu}_m))y\disc(\boldsymbol{\mu}_m), &&& z \disc \in \mathbb M, \\
& \mbf C_u: \mathbb M \rightarrow \mathbb M  \hspace{.2cm}&&
\mbf C_u(v\disc) = \frac{1}{M} \sum_{m=1}^{M} (v\disc, u\disc(\boldsymbol{\mu}_m))u\disc(\boldsymbol{\mu}_m), &&&  v \disc \in \mathbb M,\\
& \mbf C_p: \mathbb M \rightarrow \mathbb M  \hspace{.2cm}&&
\mbf C_p(q\disc) = \frac{1}{M} \sum_{m=1}^{M} (q\disc, p\disc(\boldsymbol{\mu}_m))p\disc(\boldsymbol{\mu}_m), &&& q \disc \in \mathbb M.\\
\end{align*}
Let us consider their eigenvalues $\lambda_n^y, \lambda_n^u, \lambda_n^p  \in \mbb R$ and the corresponding eigenfunctions $\xi_n^y, \xi_n^u, \xi_n^p \in \mathbb M$, with $\norm {\xi_n^y}_Y = \norm {\xi_n^u}_U = \norm {\xi_n^p}_{\red{Y}} = 1$, verifying: 
\begin{align*}
& (\mbf C_y(\xi_n^y), y\disc(\boldsymbol{\mu}_m)) = \lambda_n^y(\xi_n^y, y\disc(\boldsymbol{\mu}_m)), \hspace{1cm} 1 \leq m \leq M, \\
& (\mbf C_u(\xi_n^u), u\disc(\boldsymbol{\mu}_m)) = \lambda_n^u(\xi_n^u, u\disc(\boldsymbol{\mu}_m)), \hspace{1cm} 1 \leq m \leq M, \\
& (\mbf C_p(\xi_n^p), p\disc(\boldsymbol{\mu}_m)) = \lambda_n^p(\xi_n^p, p\disc(\boldsymbol{\mu}_m)), \hspace{1cm} 1 \leq m \leq M. 
\end{align*}
\no Let us assume that the eigenvalues satisfy $\lambda_1^i \geq \lambda_2^i \geq \cdots \geq \lambda_M^i \geq 0$, for $i=y,u,q$. The orthogonal POD basis functions are given by $\xi_1^y,\dots, \xi_M^y$, $\xi_1^u,\dots, \xi_M^u$ and $\xi_1^p,\dots, \xi_M^p$ and they span $\mathbb M$. We can take into consideration the first $N \leq M$ eigenfunctions for the sake of reduction, respectively for state, control and adjoint space: the reduced spaces $Y_N$ \red{and }  $U_N$ will be defined by them. 
%%%
%One can define the projections defined as:
%\begin{align*}
%& P_N^y : Y \rightarrow Y_N && (P_N^y[f], \xi_n^y)_Y = (f, \xi_n^y)_Y, \hspace{1cm} 1 \leq n \leq N,\\ 
%& P_N^u : U \rightarrow U_N && (P_N^u[f], \xi_n^u)_U = (f, \xi_n^u)_U, \hspace{1cm} 1 \leq n \leq N, \\
%& P_N^p: Q \rightarrow Q_N && (P_N^p[f], \xi_n^p)_Q = (f, \xi_n^p)_Q, \hspace{1cm} 1 \leq n \leq N,
%\end{align*}
%\no where $P_N^y[f], P_N^u[f]$ and $P_N^p[f]$: 
%$$
%P_N^y[f] = \sum_{n =1}^N (f, \xi_n^y)_Y \xi_n^y, \hspace{.7cm}
%P_N^u[f] = \sum_{n =1}^N (f, \xi_n^u)_U \xi_n^u, \hspace{.7cm}
%P_N^p[f] = \sum_{n =1}^N (f, \xi_n^p)_Q \xi_n^p.
%$$

%%%
\begin{remark}
\normalfont
The POD algorithm can be also seen under an algebraic point of view. For example, let us consider the control variables\footnote{The concept proposed is easily extended to state and adjoint variables, analogously.} $u\disc(\boldsymbol{\mu}_m)$ for $m = 1,\dots,M$. Let $\mbf C^u \in \mbb R^{M \times M}$ be the correlation matrix of the control \tit{snapshots}, that is:
$$
\mbf C_{mq}^u = \frac{1}{M}(u\disc(\boldsymbol{\mu}_m),u\disc(\boldsymbol{\mu}_q))_U, \hspace{1cm} 1 \leq m,q \leq M.
$$
\no Then, the $N$-largest eigenvalue-eigenvector pairs $(\lambda_n^u, v_n)$ solve the problem
$$
\mbf C^u v_n = \lambda_n v_n, \hspace{1cm} 1 \leq n \leq N, 
$$ 
\no with $\norm {v_n} = 1$. Giving a descending order to the eigenvalues $\lambda_1^u \geq \lambda_2^u \geq \cdots
\geq \lambda_N^u$, the orthogonal basis functions $\{\xi_1^u, \dots, \xi_N^u\}$ satisfy 
$U_N= \text{span }\{\xi_1^u, \dots, \xi_N^u\}$. The basis is given by:
$$
\xi_n^u = \displaystyle \frac{1}{\sqrt{M}}\sum_{m = 1}^M (v_n)_m u\disc(\boldsymbol{\mu}_m), \hspace{1cm} 1 \leq n \leq N,
$$
\no where $(v_n)_m$ is \emph{m-th} component of the control eigenvector $v_n \in \mbb R^M$.

\end{remark}
\subsection{Aggregated Reduced Spaces and Affinity Assumption}
\label{aggaff}
We now focus on the conditions needed to guarantee stability and efficiency of the proposed reduced order method. In order to prove the well-posedness of the problem \eqref{RBOCP}, the \emph{reduced inf-sup condition} of the bilinear form $\Cal B(\cdot, \cdot;\boldsymbol{\mu})$ must be fulfilled, in other words, it must exist a positive constant $\beta_{N0}$ such that
\begin{equation}
\label{infsupred}
\beta_N(\boldsymbol{\mu}) = \inf_{q_N \in \red{Y}_N} \sup_{w_N \in X_N} \frac{\Cal B(w_N,q_N; \boldsymbol{\mu})}{\norm{w_N}_X \norm{q_N}_{\red{Y}}} > \beta_{N0} \hspace{1cm} \forall \boldsymbol{\mu} \in \Cal P.
\end{equation} 

\no \red{Again, the state and the adjoint space are assumed to be the same } in order to ensure the hypothesis \eqref{infsupred}. As underlined in \cite{negri2013reduced}, simply building the reduced spaces as linear combinations of \emph{snapshots} may not lead to the fulfillment of the \emph{reduced inf-sup condition}. Then, we adopted the solution of \emph{aggregated spaces}, used in \cite{negri2015reduced,negri2013reduced}, already presented in \cite{dede2010reduced}.  This technique is based on the definition of an enriched space
$$
Z_N = \text{span }\{y\disc (\boldsymbol{\mu}^n), p\disc(\boldsymbol{\mu}^n), \; n = 1, \dots, N\}. 
$$
\no Now, let us define the reduced spaces $Y_N = Z_N$ \red{and } $X_N = Z_N \times U_N$, where
$$
U_N = \text{span }\{ u \disc (\boldsymbol{\mu} ^n), \; n = 1, \dots, N\}.
$$
\no Thanks to this choice, the \red{hypothesis of coincidence of the state and the adjoint spaces } is recovered and so the saddle-point problem \eqref{RBOCP} verifies the \emph{reduced inf-sup condition}.
\\ Let us briefly introduce the affinity assumption\footnote{If the problem does not fulfill the affinity assumption, it can be recovered thanks to the empirical interpolation method (see e.g. \cite{barrault2004empirical} and \cite[Chapter 5]{hesthaven2015certified}).} that guarantees efficiency of reduced order methods. The
problem \eqref{ZETA} admits affine decomposition if we can rewrite the bilinear forms and the functionals involved as:
\begin{equation}
\begin{matrix}
& \qquad \Cal A(x, w; \boldsymbol{\mu}) =\displaystyle  \sum_{q=1}^{Q_\Cal{A}} \Theta_\Cal{A}^q(\boldsymbol{\mu})\Cal{A}^q(x,w), & 
& \qquad \Cal B(w,p; \boldsymbol{\mu}) =\displaystyle  \sum_{q=1}^{Q_\Cal{B}} \Theta_\Cal{B}^q(\boldsymbol{\mu})\Cal{B}^q(w,p), \\
& \qquad  \la G(\boldsymbol{\mu}), s \ra =\displaystyle  \sum_{q=1}^{Q_G} \Theta_G^q(\boldsymbol{\mu})\la G^q , s \ra, & 
& \qquad  \la F(\boldsymbol{\mu}), w \ra =\displaystyle  \sum_{q=1}^{Q_F} \Theta_F^q(\boldsymbol{\mu})\la F^q , w \ra,
\end{matrix}
\end{equation}
\no  for some finite $Q_\Cal{A}, Q_\Cal{B}, Q_G, Q_F$, where $\Theta_\Cal{A}^q,\Theta_\Cal{B}^q, \Theta_G^q, \Theta_F^q$ are $\boldsymbol{\mu}-$dependent smooth functions, whereas $\Cal A^q,\Cal B^q$, \\$G^q, F^q$ are $\boldsymbol{\mu} -$independent bilinear forms and functionals. This hypothesis allows us to divide the resolution of the reduced order approximation of \eqref{ZETA} in two stages:
\begin{enumerate}
\item \textbf{offline}: in this stage the reduced spaces are built and all the $\boldsymbol{\mu}-$independent quantities are assembled. It is performed only once and it may be very costly;
\item \textbf{online}: in this phase the $\boldsymbol{\mu}-$dependent quantities are assembled and the reduced system is solved. This stage is performed every time we want the model to be simulated at a new value of $\boldsymbol{\mu}$, representing a new configuration for our system.
\end{enumerate}

\section{Applications in Environmental Marine Sciences and Engineering}
\label{app}
\begin{figure}[H]
\centering
%\subfigure[]{\includegraphics[scale = 0.25]{domains/golfomesh}}
%
\subfigure[]{\includegraphics[scale = 0.533]{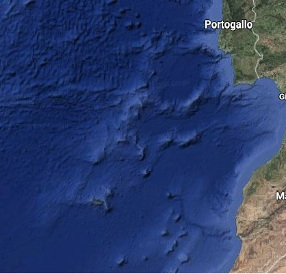}} \hspace{1cm}
\subfigure[]{\includegraphics[scale = 0.40]{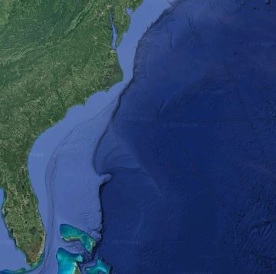}} \hspace{1cm}
\subfigure[]{\includegraphics[scale = 0.533]{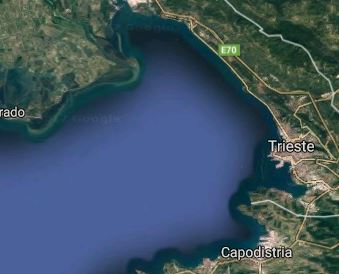}}

\caption{Satellite images. (a) North Atlantic Ocean, West coast; (b) North Atlantic Ocean, East coast; (c) Gulf of Trieste, Italy.}

\label{satellite}
\end{figure}

\begin{figure}[H]
\centering
%\subfigure[]{\includegraphics[scale = 0.25]{domains/golfomesh}}
%
\subfigure[]{\includegraphics[scale = 0.47]{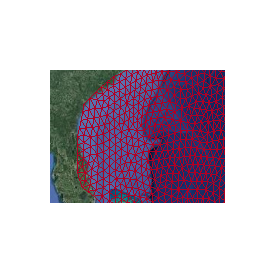}} 
\subfigure[]{\includegraphics[scale = 0.75]{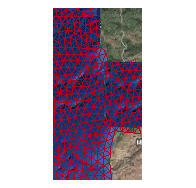}} \hspace{-1.5cm}
\subfigure[]{\includegraphics[scale = 0.47]{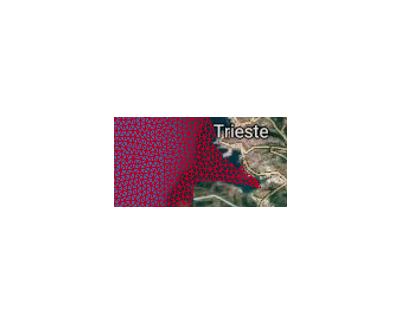}}

\caption{Meshes overlapping satellite images. (a) North Atlantic Ocean, particular of the Florida peninsula; (b) North Atlantic Ocean, particular of Portugal and North Africa; (c) Gulf of Trieste, particular of the city of Trieste.}

\label{satellitemesh}
\end{figure}

This section aims at applying the proposed reduced order method (ROM) to parametrized optimal control problems involved in environmental marine sciences and engineering problems. One of the purpose is to show the computational savings enabled by the use of a ROM in place of the usual FE approximation strategies. Two specific examples are proposed: 
\begin{enumerate}
\item \textbf{A Pollutant Control on the Gulf of Trieste}
\\ This first example involves an advection-diffusion pollutant control problem set in the Gulf of Trieste, Italy. The latter is a physical basin particularly windy and it has very peculiar \emph{flora} and \emph{fauna} population (as underlined in \cite{ mosetti2005innovative,shiganova2009native}). Moreover its analysis is important from a social point of view since it has a great impact on the local community: the city of Trieste overlooks the sea and depends on the Gulf and on its structures from harbours as well as from tourist infrastructures. For these reasons it needs to be monitored and kept under control.
\item
\textbf{A Solution Tracking of the Large Scale Ocean Circulation Model} \\ \no
The \textit{Solution tracking} is an optimal control problem that aims at making a solution the most similar to a given observation. As a second application, we propose a solution tracking problem of the large scale Ocean circulation model, governed by quasi-geostrophic equations. This OCP($\boldsymbol{\mu}$) example fits in the framework of a \textit{data assimilation technique} (see \cite{behringer1998improved, carton2000simple, carton2008reanalysis,ghil1991data} as references) that allows the model to be modified adding information from experimental data. The importance of studying Ocean Circulations Models lies in forecasting analysis of future meteorological and climatological scenarios in order to prevent catastrophic events, as described in  \cite{yang2016intensification}.
\end{enumerate} 
\no Both the presented applications are characterized by several physical parameters, and so reduced order methods could be an useful tool to decrease the time required by numerical simulations. 
\\In order to reach more realistic results, we have used meshes derived from satellite images representing the geographic area to be studied, as shown in Figure \ref{satellite}. Figure \ref{satellitemesh} gives an idea about the work needed in order to build realistic meshes: some details of the meshes overlapped to the satellite images are shown.  For our analysis, having these specific meshes was very important to give physical meaning to our problems and to their formulations, in order to achieve reliable results that could be potentially compared to real experimental data, which are strictly linked to the geographical area where they are collected. \vspace{0.5cm}
\vspace{-2cm}

\vspace{0.7cm}
\no The simulations have been run using  FEniCS \cite{fenics} for the \emph{full order} solutions and RBniCS for the reduced order ones \cite{hesthaven2015certified,rbnics}. The machine used for the simulations has a processor AMD A8-6410 APU with 8 GB of RAM.

\subsection{Reduced Basis Applied to a Pollutant Control on the Gulf of Trieste}
\begin{figure}[H]
\centering
%\subfigure[]{\includegraphics[scale = 0.25]{domains/golfomesh}}
%
\subfigure[]{\includegraphics[scale = 0.25]{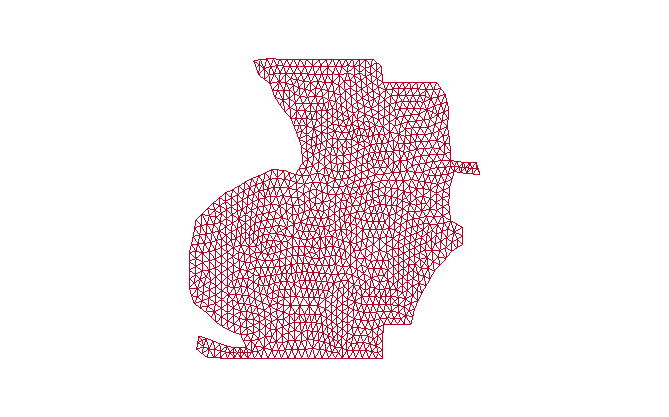}}\hspace{-1cm}
\subfigure[]{\includegraphics[scale = 0.25]{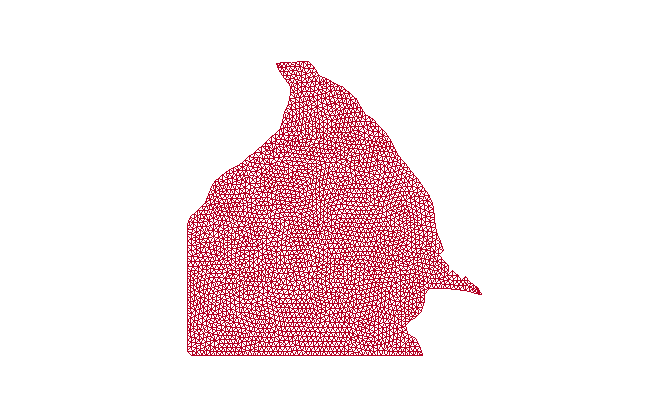}}\hspace{-1cm}
\subfigure[Boundaries: \emph{light blue}: open sea; \emph{brown}: coasts. Subdomains: \emph{red}: monitored natural area; \emph{green}: source of the pollutant tracer.]{\includegraphics[scale = 0.25]{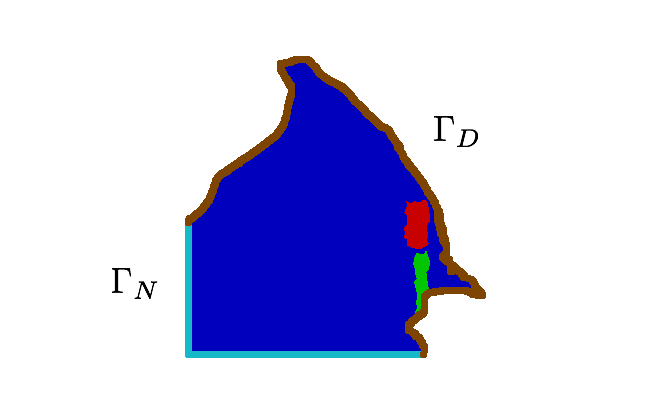}}
\caption{Mesh and subdomains. (a) mesh of the North Atlantic Ocean; (b) mesh of the Gulf of Trieste; (c) subdomains of the Gulf of Trieste.}

\label{mesh}
\end{figure}

The proposed problem aims at limiting the impact of a pollutant tracer on touristic and natural areas of the Gulf of Trieste. The OCP($\boldsymbol{\mu}$) is governed by advection-diffusion state equation (see \cite{phd, quarteroni2005numerical,quarteroni2007reduced}). Let $\Omega$ be an open, bounded and regular domain representing the Gulf of Trieste (see Figure \ref{mesh} (b)), with boundary $\partial \Omega = \Gamma_D \cup \Gamma_N$ and $\Gamma_D \cap \Gamma_N = \emptyset$, where homogeneous Dirichlet and Neumann boundary conditions are imposed on $\Gamma_D$  and $\Gamma_N$, respectively. The coasts are considered in $\Gamma_D$, while the open sea represents $\Gamma_N$ (see Figure \ref{mesh} (c)). \\ \no Let us define the state and the control spaces as $Y = H^1_{\Gamma_D}(\Omega)= \displaystyle \{ y \in H^1(\Omega) \; : 	\: y_{|\Gamma_{D}}=0\},$ and $U = \mathbb R$, respectively. \red{We remark that the adjoint space is equal to the state space}. \\ The non-dimensional OCP($\boldsymbol{\mu}$) reads: given $\boldsymbol{\mu} \in \Cal P$, find $(y(\boldsymbol{\mu}), u(\boldsymbol{\mu})) \in Y \times U$ such that:
\begin{equation}
\label{settembre}
\begin{matrix}
\underset{(y,u) \in Y \times U}{\text{min }} J(y,u) = \underset{(y,u) \in Y \times U}{\text{min }} \hspace{.2cm}
\displaystyle \half \int_{\Omega_{OBS}}(y- y_d)^2 \; d \Omega_{OBS} \;
 \red{+ \alf \int_{\Omega _u} u^2 \; d \Omega_u} \vspace{.5cm} \\
\text{such that } a(y,q; \boldsymbol{\mu}) = c(u,q), \hspace{1cm}  \forall q \in \red{Y}.
\end{matrix} 
\end{equation}

\begin{figure}[H]
%\centering
%\includegraphics[scale = 0.7]{Golfo/prova}
%\hspace{-1cm}\includegraphics[scale = 0.28]{Plots_Finali/Un_G} 
%\vspace{.7cm}\hspace{-0.5cm}\includegraphics[scale = 0.31]{Plots_Finali/FE_New_Gulf}
%\hspace{-0.5cm}\includegraphics[scale = 0.31]{Plots_Finali/RB_New_Gulf}
%\hspace{-.5cm}
\hspace{-.9cm}\includegraphics[scale = 0.25]{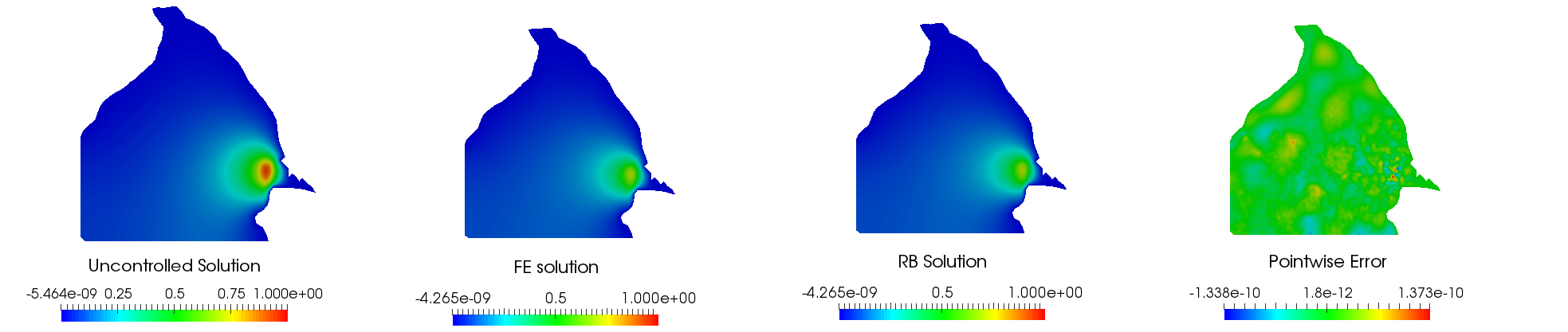}
\caption{Gulf of Trieste (Bora configuration): results.}
\label{Golfores}
\end{figure}

\begin{figure}[h]
\hspace{-1cm}\includegraphics[scale = 0.28]{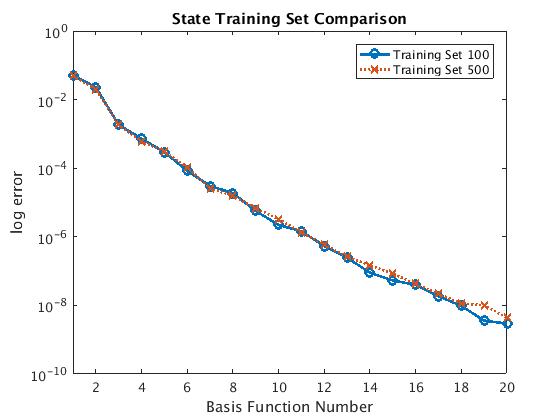} 
 \includegraphics[scale = 0.28]{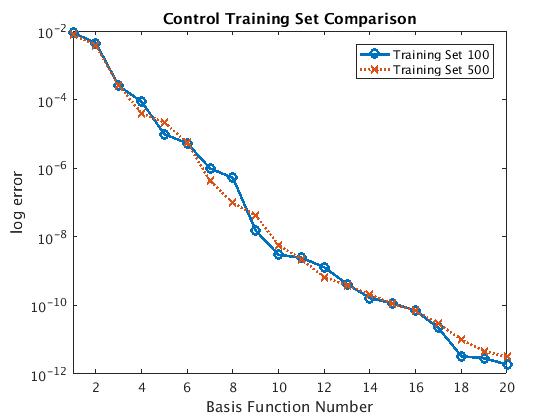} 
\includegraphics[scale = 0.28]{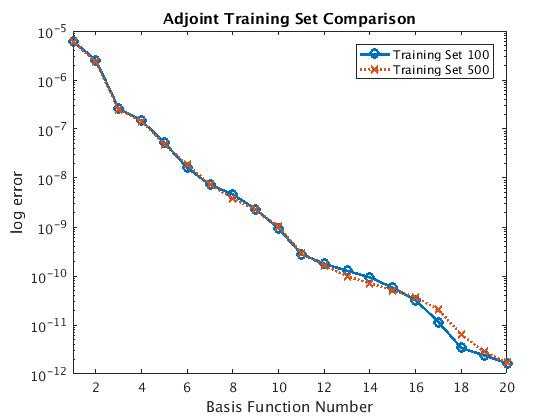} 
\caption{\reviewerA{Gulf of Trieste (Bora configuration): errors and training set comparison. The plots are almost coincident.}}
\label{golfocom}
\end{figure}

\no where the state $y$ is the pollutant concentration and $y_d = 0.2 \in \mbb R$ represents the safety threshold of the pollutant tracer. 
\no The bilinear forms $a: Y \times \red{Y} \rightarrow \mbb R$ and 
$c: U \times \red{Y} \rightarrow \mbb R$ are defined as:
\begin{align*}
\displaystyle a(y,q, \boldsymbol{\mu}) = \int_{\Omega}( \nu(\boldsymbol{\mu})\nabla y \cdot \nabla q + \boldsymbol{\beta}(\boldsymbol{\mu}) \cdot \nabla yq) \; d\Omega, &&
\displaystyle c(u,q) =  L \: u\int_{\Omega_u}q \: d\Omega_u.
\end{align*}
\no where $\nu(\boldsymbol{\mu})\equiv \mu_1$ represents the diffusivity action of the state equation, while  $\boldsymbol{\beta}(\boldsymbol{\mu})=[\beta_1(\mu_2),\beta_2( \mu_3)]$ is the transport field acting on the Gulf. Then, the parameter $\boldsymbol{\mu} = [\mu_1, \mu_2, \mu_3]$ influences the circulation of the current in the Gulf. In this case $G \in \red{Y} \dual$ is $G \equiv 0$. \red{In this experiment, the control represents the concentration of the pollutant tracer released in the green area of the Gulf (see Figure \ref{mesh} (c)) }The constant $L = 10^3$ makes the system non-dimensional. For the transport field we decided to take into consideration in proximity of the observation domain\footnote{Constant transport field will be sufficient to simulate the most interesting configurations for the circulation dynamic of the Gulf of Trieste%, and we will not take into %consideration the dynamics at the boundary
.
}$$
\beta_1(\mu_2) \equiv \mu_2, \hspace{1cm} \beta_2( \mu_3) \equiv \mu_3.
$$
The parameter space considered is $\Cal P = [0.5,1]\times[-1,1]\times[-1,1]$. The plot of the Figure \ref{mesh} (c) shows the considered subdomains: in green $\Omega_{u}$, where the pollutant loss is, in red $\Omega_{OBS}$, where we want to monitor the pollutant concentration: it represents the swimming touristic area of the city and Miramare natural area (red subdomain in Figure \ref{mesh}). Two argumentations drove us in the choice of $\Omega_{OBS}$ (as underlined in \cite{mosetti2005innovative}):
\begin{enumerate}
\item its unique ecological \emph{flora} and \emph{fauna} marine population,
\item it is an area crowded by Trieste citizens and by many tourists.
\end{enumerate}
\no The parameter $\boldsymbol{\mu}$ is a physical parameter that describes the dynamic of the currents deriving from the specific winds blowing on this geographical area (the winds acting on the Gulf will be introduced later). Varying the parameter, it is possible to simulate several configurations in order to study how the wind could affect the diffusion of a dangerous pollutant in the natural area of Miramare. 
%Reduced order methods may be a suitable tool to be used, since they allow to run many simulations, for several value of $\boldsymbol{\mu}$ in a very fast and reliable way. This is a factor of great importance, since the impact of the pollutant has to be faced as soon as possible.   
Thanks to reduced order modelling, many scenarios could be analysed, with a great saving of computational time resources. Several possible results must be taken into consideration in this context and, then, a deeper analysis can be made in order to better understand possible scenarios. The monitoring of the diffusion of a pollutant is necessary in order to preserve natural areas or to safeguard an ecological polluted habitat in case of ecological accident, and an accurate and fast model is helpful in planning a program of action.    \\ 
 \no In order to recast the problem in the framework \eqref{ZETA}, let $X = Y \times U$ be the product space of state and control spaces.  Let $x = (y,u)$ and $w = (z,v)$ be elements of $X$, whereas $q$ an element of $\red{Y}$. Moreover, we define the bilinear forms $m: Y \times Y \rightarrow \mbb R$ and $n: U \times U \rightarrow \mbb R$ as follows: 
$$
\begin{matrix}
m(y,z) = \displaystyle \int_{\Omega_{OBS}}yz \; d \Omega_{OBS} \; \;\;\text{  and } \vspace{.5cm} &
\red{n(u,v) = \displaystyle \int_{\Omega_u}uv \; d \Omega_u}. \\
\end{matrix} 
$$
Furthermore, we define the forms $\Cal A, \Cal B$ and $F$ as follows:
$$
\begin{aligned}
& \Cal A : X \times X \rightarrow \mbb R \hspace{1cm} &&  \Cal A(x,w) = m(y,z) \red{\:+ \: \alpha n(u,v)}, \\
& \Cal B : X \times \red{Y} \rightarrow \mbb R \hspace{1cm} && \Cal B(w, q; \boldsymbol{\mu})   = a(z,q; \boldsymbol{\mu}) - c(v,q),\\
& F: X \rightarrow \mathbb R \hspace{1cm} && \la F, w \ra = y_d\int_{\Omega_{OBS}}z \; d\Omega_{OBS}. \\
\end{aligned}
$$ 
To build the aggregated reduced spaces for state and adjoint we used the POD-Galerkin algorithm introduced in section 3. In this specific example the reduced space for the control does not need to be reduced and thus we set $U_N = \mathbb R$. For this problem the affinity assumption is guaranteed: with $Q_{\Cal A} = 1$, $Q_{\Cal B} =4$ and $Q_{F} =1$ the affine decomposition of the problem is given by
$$
\begin{aligned}
&\Theta^1_{\Cal A} = 1 \hspace{1cm} && \Cal A^1(x,w) = \Cal A(x,w), \\
&\Theta^1_{\Cal B} = \mu_1 \hspace{1cm} && \Cal B^1(x,q)= \displaystyle  \int_{\Omega}\nabla y \cdot \nabla q \; d \Omega,  \\
 &\Theta^2_{\Cal B} = \mu_2 \hspace{1cm} && \Cal B^2(x,q)  = \displaystyle \int_{\Omega} \frac{\partial y}{\partial x_1}q \; d\Omega,\\\
&\Theta^3_{\Cal B} = \mu_3  \hspace{1cm} && \Cal B^3(x,q)  = 
\displaystyle \int_{\Omega} \frac{\partial y}{\partial x_2}q \; d\Omega, \\
&\Theta^4_{\Cal B} = - L  \hspace{1cm} && \Cal B^4(x,q)  = 
\displaystyle \int_{\Omega_u} uq \; d\Omega, \\
&\Theta^1_F = 1 \hspace{1cm} && \la F^1, w \ra = \la F, w \ra.\\
\end{aligned}
$$
\no In the following, some numerical results linked to two different configurations are shown: the data of the experiments are reported in Table \ref{sani}. As Figure \ref{golfocom} shows, the choice of a training set of $100$ snapshots was sufficient in order to achieve a good reduced approximation: the errors\footnote{ \label{errGolfo} The errors considered for state, control and adjoint are, respectively:
$$
\norm {y\disc(\boldsymbol{\mu}) - y}_{H^1_0}, \hspace{.3cm} |u\disc(\boldsymbol{\mu}) - u_N(\boldsymbol{\mu})|, \hspace{.3cm}\norm {p\disc(\boldsymbol{\mu}) - p_N(\boldsymbol{\mu})}_{H^1_0}.
$$}  between FE and ROM variables have the same behaviour of  the solution generated by a training set of $500$ \emph{truth} approximations. We analysed how the wind action could influence pollutant diffusion. We simulated the net water transport due to Bora, a wind blowing from East to North-West, described by $(\mu_2, \mu_3) = (-1, 1)$, and Scirocco, a wind blowing from South-East, corresponding to $(\mu_2, \mu_3) = (1, -1)$. From the lower value of the cost functional reported in Table \ref {u21}, one can deduce that Bora makes the polluted water removed from $\Omega_{OBS}$, while Scirocco acts is the opposite way, pushing it towards the coast. Furthermore, Table \ref{u21} shows the differences between FE and ROM performances, in terms of dimension of the systems, time of resolution and cost functional. In Table \ref{eyein} the \emph{speed up index} behaviour with respect to the increasing of the \emph{basis numbers} is presented. The \emph{speed up index} represents the number of \emph{reduced problems} solved in the time needed for a \emph{full order} simulation. With the terminology \emph{basis number} we refer to as the number $N$ such that $N_y = N_p = 2N$ (and $N_u = 1$, being $U = \mathbb R$), where $N_y, N_u, N_p$ are the number of online degrees of freedom for state, control and adjoint variable, respectively. In other words if the \emph{basis number} is $N$, we are solving a system of dimension $(4N +1)  \times (4N + 1)$.
\\
\no The first left plot in Figure \ref{Golfores}  shows an ``uncontrolled'' Bora configuration (corresponding to $\boldsymbol{\mu} = (1, -1, 1)$, resulting from the simulation of the state equation (only) with a value $u=1$ as a forcing term. Then, in the same Figure \ref{Golfores}, the optimal control problem solution is presented with FE discretization and ROM, respectively. As one can see, the FE approximation and the ROM one match. Another proof of the reliability of reduced basis POD-Galerkin method is the pointwise error shown in the last plot of Figure \ref{Golfores}: the maximum value reached is $1.373 \cdot 10^{-10}$ with \emph{basis number} $N = 20$.

%\begin{minipage}{0.9\textwidth}
%\end{minipage}

%\footnote{\label{speed}The speed up index represents the number of reduced problems solved in the time needed for a FE simulation.} 
\vspace{-0.5cm}
\begin{table}[h]
\centering
\caption{Data of the numerical experiment: Gulf of Trieste.}
\label{sani}
\begin{tabular}{ c | c | c}
\toprule
\textbf{Data} & \textbf{Bora Values} & \textbf{Scirocco Values} \\
\midrule
$\mu_1$ & 1& 1 \\

 $(\mu_2, \mu_3)$ &  $( -1, 1)$ & $(1, -1)$ \\
\midrule
$y_d$  & 0.2 & 0.2\\ \midrule
POD Training Set Dimension & 100 & 100\\ 
\midrule
Basis Number $N$ & 20 & 20\\
\midrule
Sampling Distribution & uniform & uniform \\
\bottomrule
\end{tabular}
\vspace{-.1cm}
\begin{table}[H]
\centering
\caption{ROM vs FE:  Gulf of Trieste (Bora configuration).}
\label{u21}
\begin{tabular}{c | c| cl}
\toprule
\textbf{Bora Configuration} & FE & ROM \\
\midrule
System Dimension &   $5939 \times 5939$  & $201 \times 201$ \\
Optimal Cost Functional & $4.9167 \cdot 10^{-5}$ & $4.9167 \cdot 10^{-5}$ \\
Time of Resolution & $2.79s$ & $2.41 \cdot 10^{-2}s$ \\
\midrule
\textbf{Scirocco Configuration} & FE & ROM \\
\midrule
System Dimension & $5939 \times 5939$  &   $201 \times 201$ \\
Optimal Cost Functional & $5.3417 \cdot 10^{-5}$ & $5.3417 \cdot 10^{-5}$ \\
Time of Resolution & $3.12s$ & $3.41 \cdot 10^{-2}s$  \\
\bottomrule
\end{tabular}
\end{table}
\end{table}
\vspace{-.5cm}
\begin{table}[H]
\centering
\caption{Speed up analysis: Gulf of Trieste (Bora configuration).}
\label{eyein}
\begin{tabular}{c | cccccccccccl}
\toprule

Basis Number  $N$            & 1 & 5 & 10 & 5 &20 \\
\reviewerA{Speed up}     &         \reviewerA{361}& \reviewerA{364} & \reviewerA{350} &\reviewerA{ 317} & \reviewerA{296} \\
%$\Cal N \times \Cal N$   & $1323 \times 1323$ &   &   &   &   &   &   &   &   &       \\
%$N \times N$ & $100 \times 100$   &   &   &   &   &   &   &   &   &     \\
%Basis Number  $N$            & 1                & 2 & 3 & 4 & 5 & 6 & 7 & 8 & 9 & 10  \\%& 11 & 12 & 13 & 14 & 15 & 16 & 17 & 18 & 19 & 20 \\
%Speed up     &                   
%361&
%380&
%369&
%366&
%364&
%362&
%352&
%348&
%346&
%350\\
%\midrule
%Basis Number $N$  & 11 & 12 & 13 & 14 & 15 & 16 & 17 & 18 & 19 & 20 \\
%Speed up                              &    
%336&
%334&
%333&
%328&
%317&
%315&
%309&
%303&
%301&
%296\\
%\midrule
%Basis Number $N$ & 21 & 22 & 23 & 24 & 25 & 26 & 27 & 28 & 29 & 30 \\
%%%& 31 & 32 & 33 & 34 & 35 & 36 & 37 & 38 & 39 & 40 &
%%%& 41 & 42 & 43 & 44 & 45 & 46 & 47 & 48 & 49 & 50 &
%Speed up                              &    
%294&
%288&
%282&
%277&
%271&
%264&
%267&
%257&
%251&
%245\\
%\midrule
%Basis Number $N$ &31 & 32 & 33 & 34 & 35 & 36 & 37 & 38 & 39 & 40\\
%%%& 31 & 32 & 33 & 34 & 35 & 36 & 37 & 38 & 39 & 40 &
%%%& 41 & 42 & 43 & 44 & 45 & 46 & 47 & 48 & 49 & 50 &
%Speed up                              &    
%240&
%234&
%231&
%218&
%216&
%204&
%201&
%202&
%195&
%188\\
%\midrule
%Basis Number $N$ & 41 & 42 & 43 & 44 & 45 & 46 & 47 & 48 & 49 & 50\\
%%%& 31 & 32 & 33 & 34 & 35 & 36 & 37 & 38 & 39 & 40 &
%%%& 41 & 42 & 43 & 44 & 45 & 46 & 47 & 48 & 49 & 50 &
%Speed up                              &    
%183&
%175&
%167&
%161&
%157&
%152&
%147&
%147&
%134&
%135\\
\bottomrule
\end{tabular}
%\\
%\vspace{.5cm}
%In the following a comparison between the dimensions of full order system and the reduced one is shown.
%\vspace{.5cm}\\
%\begin{tabular}{c | c}
%\toprule
%$\Cal N \times \Cal N$   & $5639 \times 5639$       \\
%\midrule
%$4N + 1 \times 4N + 1$ & $201 \times 201$   \\
%\bottomrule
%\end{tabular}
\end{table}

\begin{figure}[H]
\centering
\includegraphics[scale = 0.31]{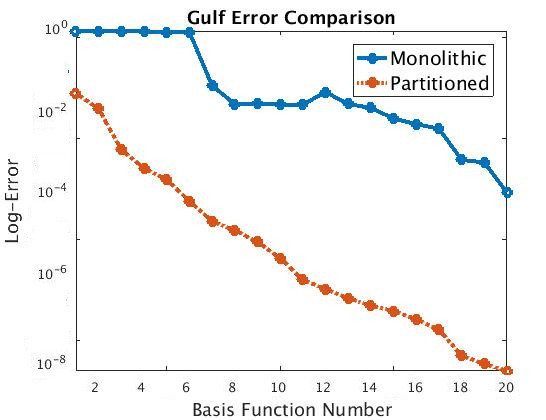}

\caption{\reviewerB{Gulf of Trieste (Bora configuration): errors.}}
\label{golfoerr}
\end{figure}
 \no In Figure \ref{golfocom} errors between FE solution and ROM solution with respect to the \emph{basis number} $N$ over a random testing set of $100$ are presented. The obtained results show that the ROM allows to get fast and accurate simulations, since very few basis functions are required to have very small errors.
\\The plot in Figure \ref{golfoerr} illustrates a comparison between the POD approach presented in section 3 (labeled as \emph{partitioned}) and a different reduction in which only one POD is carried out for the variables $(y, u, p)$ at the same time on the space $X \times \red{Y}$ (labeled as \emph{monolithic}).
The results show that it is preferable to use a partitioned POD approach rather than a monolithic POD algorithm. The improvement is significant, since the latter approach gave an error of the order of $10^{-4}$ with $N = 20$, while, for the same value of $N$ and with the partitioned option, the sum of the state, control and adjoint errors reaches the value of $10^{-8}$. Since the partitioned approach gave better values of the errors, we decided to exploit it for the Oceanographic application which we will present in the next section.
\subsection{Reduced Basis Applied to an Ocean Circulation Solution Tracking}
\label{OCEANO}
\begin{figure}[h]

\hspace{-1.5cm}\includegraphics[scale = 0.21]{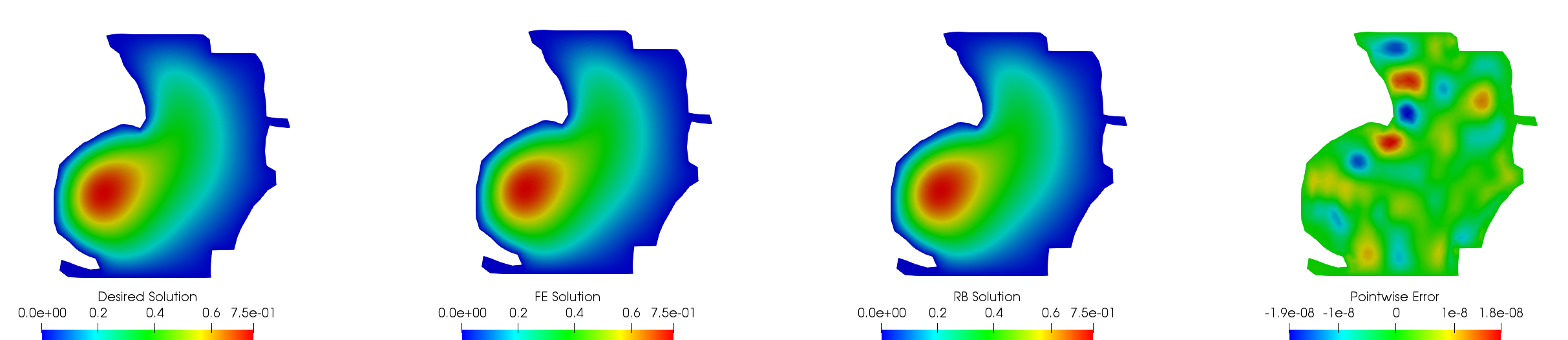}
%\hspace{-1.2cm}\includegraphics[scale = 0.24]{Atlantico/FE_solution}
%\hspace{-1.6cm}\includegraphics[scale = 0.24]{Atlantico/RB_solution}\hspace{-1.7cm}\includegraphics[scale = 0.24]{Atlantico/DIFF}
\caption{North Atlantic Ocean linear solution tracking: results.}
\label{NAres}
\end{figure}

\begin{figure}[H]

\centering
\hspace{-.8cm}
\includegraphics[scale = 0.28]{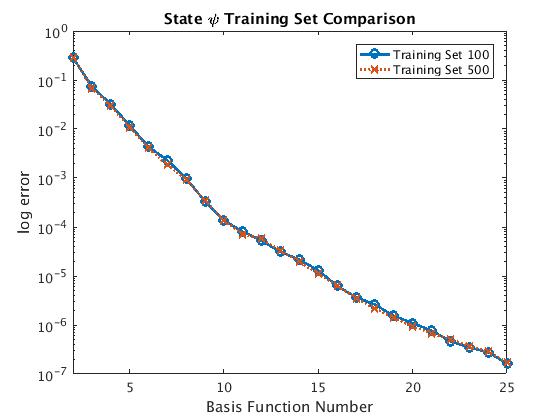}\includegraphics[scale = 0.28]{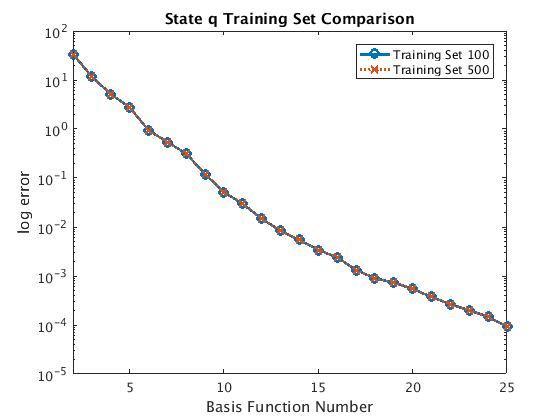}\includegraphics[scale = 0.28]{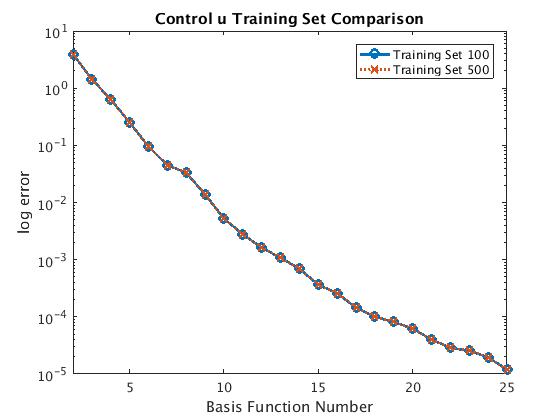}\\ 
\includegraphics[scale = 0.28]{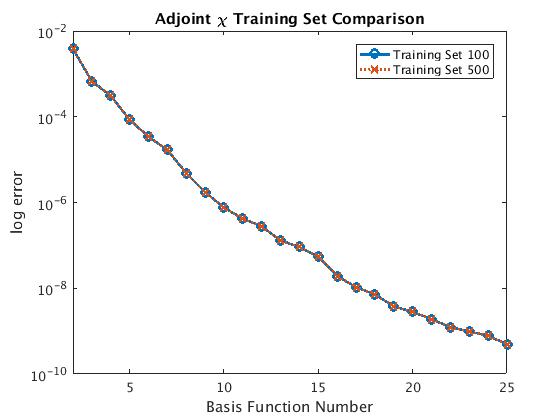}\includegraphics[scale = 0.28]{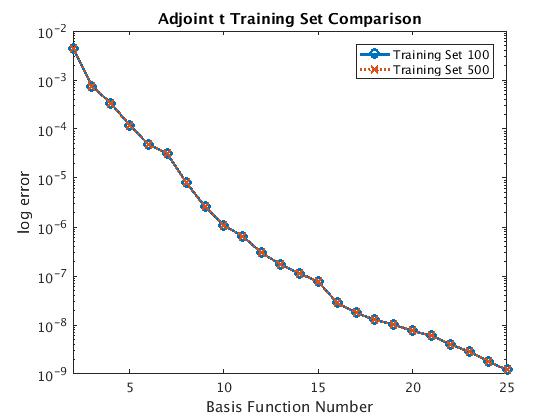}

\caption{\reviewerA{North Atlantic Ocean linear solution tracking: errors and training set comparison. The plots are coincident.}}
\label{NAcom}
\end{figure}

The general Ocean circulation model describes large scale flow dynamics. It is strictly linked to wind action: it can be considered as a coupled system of Ocean and Atmosphere. The theory associated to this topic is deeply analysed in \cite[Chapter 3]{cavallini2012quasi}. The model is governed by the following non-dimensional PDE($\boldsymbol{\mu}$), known as \textit{quasi-geostrophic} equation:
\begin{equation} 
\label{QG} 
\mu_3 \mathcal{F}(\psi, \Delta\psi)+ \frac{\partial \psi}{\partial x} = f  - \mu_1 \Delta\psi + \mu_2 \Delta^2\psi,
\end{equation}
\no where, given a suitable spaces $V$, the nonlinearity of the expression is given by $\Cal F: V \times V \rightarrow \mathbb R$ defined as:
\begin{equation}
\label{QGnl}
\Cal F (\psi, q) = \frac{\partial \psi}{\partial x}\frac{\partial q}{\partial y} 
- \frac{\partial \psi}{\partial y}\frac{\partial q}{\partial x}.
\end{equation}
The parameter $\mu_3$ represents how much the nonlinear term affects the flow dynamics, while $\mu_1$ and $\mu_2$ stand for diffusive action, respectively.
 \\The parameter $\boldsymbol{\mu}$ describes the North Atlantic Ocean dynamics completely, since it gives information about how the large scale Ocean circulation is affected by different phenomena, such as location (typically described by $\mu_3$) and intensity variations of the gyres and the currents in the Ocean. Let us recall that Ocean dynamic is strongly linked to wind stress and atmospheric behaviour. Then, the parameter $\boldsymbol{\mu}$ describes the dynamic of a very complex physical system, taking into account several natural factors and phenomena. It is very important to run many simulations for different values of the parameter $\boldsymbol{\mu}$, in order to achieve a full knowledge of this system describing Oceanic and Atmospheric dynamics, linked to climatological forecasting issues.  \\ \no As specified in \cite[Section 3.2]{cavallini2012quasi}, the forcing term $f$ depends on wind stress $\boldsymbol{\tau}$  by the following relation:
$$
f = \hat{\mbf k} \cdot \text{rot } {\boldsymbol{\tau}},
$$
\no where $\hat{\mbf k}$ is the third reference spatial unit vector. In our application we considered a bounded and regular bi-dimensional domain\footnote{Experiments showed that the analysis of the underwater depth did not affect the dynamics of the equation, so we decided to exploit a simpler bi-dimensional model. Although, adding the bathimetric effect is quite simple: let us suppose that the ocean floor is described by a smooth function $h: \Omega \rightarrow \mathbb R$, then one can consider
$$  \mu_3 \mathcal{F}(\psi, \Delta\psi + h)+ \frac{\partial \psi}{\partial x} = f  - \mu_1 \Delta\psi + \mu_2 \Delta^2\psi,$$ 
in order to treat a more complete model. As one can see, the bathimetry affects only the nonlinear term.} $\Omega \subset \mathbb R^2$ representing the North Atlantic Ocean (see Figure \ref{mesh} (a)). 
 \no Furthermore, \nl{in this section, } we focused on the linear version of the quasi-geostrophic equation ($\mu_3 = 0$). \nl{The nonlinear $OCP(\pmb \mu)$ governed the oceanographic model will be treated in section \ref{nl_sec}. }The solution tracking problem constrained to this particular state equation is:

\begin{equation}
\label{nodec}
\begin{aligned}
& \underset{(\psi, u) \in V \times U}{\text{min }}
 J(\psi, u) =\underset{(\psi, u) \in V \times U}{\text{min }} \hspace{.2cm} \displaystyle \half \int_{\Omega} (\psi - \psi_d)^2 \; d \Omega + \displaystyle \alf \int_{\Omega} u^2 \; d \Omega \vspace{0.2cm}\\
& \qquad \text{such that } \hspace{0.5cm}
\begin{cases}
\displaystyle \frac{\partial \psi}{\partial x} = u
- \displaystyle \mu_1 \Delta \psi + \mu_2 \Delta ^2 \psi & \text{in } \Omega, \\
\psi = 0 & \text{on } \partial \Omega, \\
\Delta \psi = 0 & \text{on } \partial \Omega,
\end{cases}
\end{aligned}
\end{equation}
\no where $\psi \in V$ is our state variable, $u \in U$ is the forcing term to be controlled, where $V$ and $U$ are two suitable functions spaces and $\alpha$ is the penalization term. The physical parameter $\boldsymbol{\mu} = (\mu_1, \mu_2)$ is considered in the parametrized space $\Cal P = [10^{-4},1]\times [10^{-4},1]$.  \red{In this example, the control variable represents the wind action. We stress that the quasi--geostrophic model describes how Ocean and Atmosphere interact. }The problem aims at making the solution $\psi$ the most similar to $\psi_d$, \red{representing the desired Ocean dynamics, thanks to the action of to the wind stress expressed by $u$, as mentioned above}.
\\ In some applications $\psi_d$ can \reviewerA{represent} experimental data. In this sense, our experiment could be seen as a prototype of a \textit{data assimilation} model with forecasting purposes (see \cite{behringer1998improved,carton2000simple, carton2008reanalysis, ghil1991data, kalnay2003atmospheric, tziperman1989optimal} as references).  This particular technique changes the model in order to achieve a solution comparable with real experimental data. \textit{Data assimilation} techniques are very costly, and reduced order methods fit very well in this context. The proposed experiment helps to understand how reduced order techniques could be exploited in order to simulate several climatological scenarios in a low dimensional and accurate framework. The opportunity of running the reduced model many times allow us to have a deeper comprehension of the Ocean dynamic, and of atmospheric phenomena and climatological scenarios. 
\\ \no In order to manage a handier problem\footnote{ \label{coe}This version of the problem does not ensure the coercivity of the state equation, that is proved in \cite{kim2015b} for the state equation of the system \eqref{nodec}}, we rewrite the previous system as it follows: 
\begin{equation}
\label{dec}
\begin{aligned}
& \underset{((\psi ,q), u) \in Y \times U}{\text{min }}
 J((\psi, q), u) =\underset{((\psi ,q), u) \in Y \times U}{\text{min }}\hspace{.2cm} \displaystyle \half \int_{\Omega} (\psi - \psi_d)^2 \; d \Omega + \displaystyle \alf \int_{\Omega} u^2 \; d \Omega \vspace{0.2cm}\\
& \qquad \text{such that } \hspace{0.5cm}
\begin{cases}
q = \Delta \psi  & \text{in } \Omega, \\
\displaystyle \frac{\partial \psi}{\partial x} = +u
- \displaystyle \mu_1 q + \mu_2 \Delta q & \text{in } \Omega, \\
\psi = 0 & \text{on } \partial \Omega, \\
q = 0 & \text{on } \partial \Omega,
\end{cases}
\end{aligned}
\end{equation}
\no where the spaces are defined as $Y = H^1_0(\Omega)\times H^1_0(\Omega)$ and $U = L^2(\Omega)$. The weak formulation reads as: 
\begin{equation}
\label{statodebole}
a((\psi, q), (\phi, r); \boldsymbol{\mu}) = c(u,(\phi, r)) \hspace{1cm} \forall \phi, r \in H^1_0(\Omega),
\end{equation}
\no where $a : Y \times Y \rightarrow \mbb R$ and $c: U \times Y \rightarrow \mbb R$ are given by:
\begin{align}
\label{stato1}
\displaystyle a((\psi, q), (\phi, r); \boldsymbol{\mu}) & = \int_{\Omega} \frac {\partial \psi }{\partial x} r \; d \Omega
\displaystyle + \mu_2  \int_{\Omega} \nabla q \cdot \nabla r \; d \Omega \;
 \displaystyle +\mu_1 \int_{\Omega} q r \; d \Omega 
 \displaystyle + \int_{\Omega} q \phi \; d \Omega
\displaystyle + \int_{\Omega} \nabla \psi \cdot \nabla \phi \; d \Omega, \vspace{0.2cm} \\ \label{stato2} %\tag{4.2.4} \\
\displaystyle c(u,(\phi, r)) & = \int_{\Omega} u r \; d \Omega .
\end{align}
\no In this case $G \in \red{Y}\dual$ is $G \equiv 0$. \vspace{.5cm} \\
Since we are facing a linear quadratic optimal control problem, it can be recast in saddle-point formulation \eqref{ZETA}. Let us define the product space 
$X = Y \times U$ and let $x = ((\psi, q),u)$ and $w = ((\chi, t),v)$ be two elements of $X$, whereas let $s = (\phi, r)$ be an element of the adjoint space. Furthermore, one has to specify:
$$
\begin{aligned}
& \Cal A : X \times X \rightarrow \mbb R \hspace{1cm} &&  \Cal A(x,w) = m((\psi, q), (\chi,t)) + \alpha n(u,v),  \\
& \Cal B : X \times \red{Y} \rightarrow \mbb R \hspace{1cm} && \Cal B(w, s; \boldsymbol{\mu})  = a((\chi,t),(\phi, r), \boldsymbol{\mu}) - c(v, (\phi,r)), \\
& F: X \rightarrow \mathbb R \hspace{1cm} && \la F, w \ra = \int_{\Omega} \psi_d \chi \; d \Omega, \\
\end{aligned}
$$
where $m : Y \times Y \rightarrow \mbb R$ and 
$n: U \times U \rightarrow \mbb R$ are defined as
$$
\begin{matrix}
\displaystyle m((\psi, q), (\chi,t)) = \int_{\Omega} \psi \chi \; d \Omega, &
\displaystyle n(u,v) = \int_{\Omega}uv \; d \Omega.\\
\end{matrix}
$$
\no In order to build the aggregated reduced spaces of the type proposed in section \ref{aggaff} a POD-Galerkin algorithm  has been exploited. 
As we underlined in section \ref{aggaff}, the affinity assumption must be guaranteed for the efficiency of the reduced problem. Indeed, with $Q_{\Cal A} = 1$, $Q_{\Cal B} =\reviewerA{3}$ and $Q_{F} =1$ the affine decomposition of the problem is given by
$$
\begin{aligned}
&\Theta^1_{\Cal A} = 1 \hspace{1cm} && \Cal A^1(x,w) = \Cal A(x,w), \\
&\Theta^1_{\Cal B} = \mu_1 \hspace{1cm} && \Cal B^1(x,s)= \displaystyle  \int_{\Omega} q r \; d \Omega,  \\
 &\Theta^2_{\Cal B} = \mu_2 \hspace{1cm} && \Cal B^2(x,s)  = \displaystyle \int_{\Omega} \nabla q \cdot \nabla r \; d \Omega,\\
 &\Theta^3_{\Cal B} = 1 \hspace{1cm} && \Cal B^3(x,s)  = \displaystyle\int_{\Omega} \frac {\partial \psi }{\partial x} r \; d \Omega
+ \displaystyle  \int_{\Omega} q \phi \; d \Omega
+ \int_{\Omega} \nabla \psi \cdot \nabla \phi \; d \Omega - \int_{\Omega}u r \; d \Omega,\\
&\Theta^1_F  = 1 \hspace{1cm} && \la F^1, w \ra = \la F, w \ra.\\
\end{aligned}
$$
%\no \textcolor{red}{\textbf{inizio aggiunta}}\\
%Since we are facing an optimal control problem governed by an elliptic coercive state equation, the well-posedness of the problem can be proved, in other words, conditions (i)-(iv) presented in section \ref{LQ} are fulfilled, as specified in \cite{negri2013reduced}. We focus our attention on the properties of the bilinear form $a(\cdot, \cdot)$. The properties needed by the forms $c(\cdot, \cdot)$, $n(\cdot, \cdot)$ and $n(\cdot, \cdot)$ to fulfill (i)-(iv) (they are presented in \cite{bochev2009least}) are trivial. \textcolor{blue}{\textbf{CREDETE CHE DOVREI ELENCARLE??}}
%Our purpose is to verify the hypotheses (i)-(iv) introduced in section \ref{conty}, in order to prove the well-posedness of the solution tracking optimal control problem. Let us %recall that inequality (D) holds for $X$ and $Q$. The following holds:
%\begin{enumerate}[(i)]
%\item the bilinear form  $\Cal A(\cdot, \cdot; \boldsymbol{\mu})$ is continuous over $X \times X$, since
%\begin{align*}
%|\Cal {A}(x, w)| & \leq \norm{\psi}_2\norm{\chi}_2 + \alpha \norm{u}_2\norm{v} \\
%			& \leq \max \{1,\alpha \} \norm{x}_X \norm{w}_X; \\
%\end{align*}
%\item
%\item the bilinear form  $\Cal B(\cdot, \cdot; \boldsymbol{\mu})$ is continuous over $X \times Q$, since
%\begin{align*}
%|\Cal {B}(x, s; \boldsymbol{\mu})| & \leq | a((\psi, q), (\phi, r); \boldsymbol{\mu})|+ |c(u, (\phi, r))| \\
%					& \leq \max \{ \mu_1, \mu_2, 1\}\norm{x}_X\norm{s}_Q;
%\end{align*}
%\item 
%\end{enumerate}
% \textcolor{red}{\textbf{fine aggiunta}}\\
\begin{table}[H]
\centering
\caption{Data of the numerical experiment: North Atlantic Ocean linear solution tracking.}
\label{antony}
\begin{tabular}{ c | c }
\toprule
\textbf{Data} & \textbf{Values} \\
\midrule
 $(\mu_1, \mu_2, \alpha)$ &  $(10^{-4}, 0.07^3, 10^{-5})$ \\
\midrule
& FE solution of \\ $\psi_d$  & quasi-geostrophic equation \\ & with $f = -\sin (\pi y)$\\
& and $\pmb \mu = (10^{-4}, 0.07^3)$\\ \midrule
POD Training Set Dimension & 100 \\ 
\midrule
Basis Number $N$ & \reviewerA{25}\\
\midrule
Sampling Distribution & log-uniform \\
\bottomrule
\end{tabular}
\end{table}
\no In Table \ref{antony} the data of the experiment are shown: \reviewerA{ choosing a training set of 100 generating elements gives comparable results with respect the one achieved with a training set generated by 500 snapshots, as presented in the errors\relax\footnote{\label{error}\reviewerA{The errors considered for state, control and adjoint are, respectively:
$$
\norm {\psi\disc(\boldsymbol{\mu}) - \psi_N(\boldsymbol{\mu})}_{H^1_0}, \hspace{.3cm} \norm {q\disc(\boldsymbol{\mu}) - q_N(\boldsymbol{\mu})}_{H^1_0}, \hspace{.3cm} \norm {u\disc(\boldsymbol{\mu}) - u_N(\boldsymbol{\mu})}_{L^2}, \hspace{.3cm}\norm {\chi\disc(\boldsymbol{\mu}) - \chi_N(\boldsymbol{\mu})}_{H^1_0},\hspace{.3cm}\norm {t\disc(\boldsymbol{\mu}) - t_N(\boldsymbol{\mu})}_{H^1_0}.
$$}} plotted in Figure \ref{NAcom}}. In Figure \ref{NAres} the desired $\psi_d$ value to be reached is presented with the FE and ROM solutions. The approximated solutions match. The last plot of the Figure \ref{NAres} shows the pointwise error: the maximum value reached is $1.8 \cdot 10^{-8}$ with \emph{basis number} $N = 25$.

%\begin{figure}[H]
%\centering
%\includegraphics[scale = 0.3]{Atlantico/D_solution}\includegraphics[scale = 0.3]{Atlantico/FE_solution}
%\\ \includegraphics[scale = 0.3]{Atlantico/RB_solution}\includegraphics[scale = 0.3]{Atlantico/DIFF}
%\caption{North Atlantic Ocean solution tracking: results.}
%\label{NAres}
%\end{figure}
%\vspace{-.8cm}
\no The following tables represent the comparison between ROM and FE performances, in terms of: system dimension, cost functional optimal value, time of resolution (Table \ref{u2}) and \textit{speed up} index with respect to the \textit{basis number} $N$ such that $2 N = N_{\psi} = N_{q} = N_{\chi} = N_{t}$ and $N = N_u$ (Table \ref{signoraluna}), which results in a solution of a $9N \times 9N$ linear system. We can deduce how \reviewerA{the ROM} method is a suitable and convenient approach to study large scale phenomena in oceanography, a field dealing with parametrized simulations that require days of CPU times for complex configurations.

\begin{table}[h]
\centering
\caption{ROM vs FE: North Atlantic Ocean linear solution tracking.}
\label{u2}
\begin{tabular}{c | c| cl}
\toprule
& FE & ROM \\
\midrule
System Dimension & $6490\times 6490$  & $450 \times 450$ \\
Optimal Cost Functional & $1.520 \cdot 10^{-6}$ & $1.520 \cdot 10^{-6}$ \\
Time of Resolution & $6.07s$ & $2.03 \cdot 10^{-1}s$ \\
\bottomrule
\end{tabular}
\end{table}
\vspace{-.4cm}
\begin{table}[H]
\centering
\caption{Speed up analysis: North Atlantic Ocean linear solution tracking.}
\label{signoraluna}
\begin{tabular}{c | cccccccccccl}

\toprule

\reviewerA{Basis Number } $\reviewerA{N}$            &\reviewerA{ 1} & \reviewerA{5} & \reviewerA{10} & \reviewerA{15} &\reviewerA{20}& \reviewerA{25} \\
\reviewerA{Speed up}     &       \reviewerA{  1049}& \reviewerA{395} & \reviewerA{192} &\reviewerA{207} & \reviewerA{199} & \reviewerA{33} \\
%%$\Cal N \times \Cal N$   & $1323 \times 1323$ &   &   &   &   &   &   &   &   &       \\
%%$N \times N$ & $100 \times 100$   &   &   &   &   &   &   &   &   &     \\
%Basis Number $N$          & 1                & 2 & 3 & 4 & 5 & 6 & 7 & 8 & 9 & 10  \\%& 11 & 12 & 13 & 14 & 15 & 16 & 17 & 18 & 19 & 20 \\
%Speed up     &                   
%1049&
%369&
%437&
%309&
%395&
%223&
%315&
%385&
%247&
%192\\
%\midrule
%Basis Number $N$ & 11 & 12 & 13 & 14 & 15 & 16 & 17 & 18 & 19 & 20 \\
%Speed up                              &    
%263&
%342&
%253&
%253&
%207&
%210&
%205&
%193&
%180&
%199\\
%\midrule
%Basis Number $N$ & 21 & 22 & 23 & 24 & 25 & 26 & 27 & 28 & 29 & 30 \\
%%%& 31 & 32 & 33 & 34 & 35 & 36 & 37 & 38 & 39 & 40 &
%%%& 41 & 42 & 43 & 44 & 45 & 46 & 47 & 48 & 49 & 50 &
%Speed up                              &    
%132&
%151&
%95&
%66&
%33&
%75&
%90&
%107&
%66&
%54\\
%\midrule
%Basis Number $N$ &31 & 32 & 33 & 34 & 35 & 36 & 37 & 38 & 39 & 40\\
%%%& 31 & 32 & 33 & 34 & 35 & 36 & 37 & 38 & 39 & 40 &
%%%& 41 & 42 & 43 & 44 & 45 & 46 & 47 & 48 & 49 & 50 &
%Speed up                              &    
%48&
%54&
%80&
%78&
%62&
%48&
%53&
%52&
%51&
%44\\
%\midrule
%Basis Number $N$ & 41 & 42 & 43 & 44 & 45 & 46 & 47 & 48 & 49 & 50\\
%%%& 31 & 32 & 33 & 34 & 35 & 36 & 37 & 38 & 39 & 40 &
%%%& 41 & 42 & 43 & 44 & 45 & 46 & 47 & 48 & 49 & 50 &
%Speed up                              &    
%29&
%41&
%35&
%34&
%31&
%28&
%25&
%21&
%17&
%23\\
\bottomrule
\end{tabular}
%\vspace{.5cm}
%\\ In the following a comparison between the dimensions of full order system and the reduced one is shown.
%\vspace{.5cm}\\
%\begin{tabular}{c | c}
%\toprule
%$\Cal N \times \Cal N$   & $5935 \times 5935$       \\
%\midrule
%$9N \times 9N$ & $450 \times 450$   \\
%\bottomrule
%\end{tabular}

\end{table}

%\begin{figure}[H]
%\centering
%\hspace{-.5cm}
%\includegraphics[scale = 0.28]{PlotOcean/psi}\includegraphics[scale = 0.28]{PlotOcean/q}\includegraphics[scale = 0.28]{PlotOcean//u}\\ 
%\includegraphics[scale = 0.28]{PlotOcean/chi}\includegraphics[scale = 0.28]{PlotOcean/t}
%
%\caption{North Atlantic Ocean linear solution tracking: errors.}
%\label{NAerr}
%\end{figure}
%\vspace{-.5cm}
\no In Figure \ref{NAcom} the error norm between the FE approximation and the ROM discretization over a random testing set of $100$ is shown for all the variables, state $\psi$ and $q$, control $u$, and adjoint $\chi$ and $t$, respectively.
\nl{
\section{Nonlinear version of the Ocean Circulation Solution Tracking}
\label{nlver}

This section aims at introducing a nonlinear version of the Oceanographic solution tracking proposed in section \ref{OCEANO}. First of all, the problem will be presented in a general theoretical formulation, then the $OCP(\pmb \mu)$ will be specified to the nonlinear quasi-geostrophic equations. As we did for the linear case, some numerical results will be presented, in order to sustain the idea of the great versatility of the ROMs in these kind of applications. 

}
\nl{
\subsection{Nonlinear reduced OCP($\pmb \mu$)s: a brief introduction}

In the following section we want to briefly describe nonlinear OCP($\pmb \mu$)s, in order to make the reader acquainted with the nonlinear numerical example presented. Let us  introduce all the quantities needed, following the structure already exploited in section
\ref{LQ}. We remark that, as in the linear case, in the different definitions, the dependence from the parameter $\pmb{\mu}$ is sometimes understood.
\\ Let $\Omega \subset \mathbb{R}^n$ be an open regular domain of boundary $\partial \Omega$. Let us introduce the Hilbert spaces $Y$, $U$ and $Z \supseteq Y$, where $Y$ is used both for the state and for the adjoint space, $U$ is the control space and $Z$ is the observation space. With a totally analogous procedure with respect the one exploited in \eqref{eq : stato}, we can define the nonlinear constraint equation: 
\begin{equation}
\label{eq : statonl}
a_{\text{nl}}(y,q; \boldsymbol{\mu}) = c(u,q; \boldsymbol{\mu}) + \la G(\boldsymbol{\mu}), q \ra \hspace{1cm} \forall q \in Y,
\end{equation}

\no where $a_{\text{nl}} : Y \times Y \rightarrow \mathbb R$ represents a continuous nonlinear state operator, $c : U \times Y \rightarrow \mathbb R$ is the usual continuous bilinear form linked to the control variable and $G \in Y\dual$ gives us information about forcing and boundary terms. The subscript ``nl'' highlights the nonlinearity of the governing state equation.
Set a constant $\alpha > 0$, we can consider again the linear quadratic objective functional \eqref{obj}
where the roles of $y_d(\boldsymbol{\mu}) \in Z$, $m : Z \times Z \rightarrow \mathbb R$ and $n : U \times U \rightarrow \mathbb R$ remain the same of the linear case in section \ref{LQ}. The OCP($\pmb \mu$) reads as the problem \eqref{eq: controllo}, with the only difference that the constraint is nonlinear.\\
Let us define the product space $X = Y \times U$, endowed with the same scalar product and the same norm introduced in section \ref{LQ}. Taking $x = (y,u)$ and $w = (z,v)$ as two elements of $X$, we can define the following quantities:
\begin{align*} & \Cal A : X \times X \rightarrow  \mathbb R &&
\Cal A (x, w; \boldsymbol{\mu}) = m(y,z; \boldsymbol{\mu}) + \alpha n(u,v; \boldsymbol{\mu}) \hspace{1cm} &&& \forall x , w \in X, \\
& \Cal B_{\text{nl}} : X \times Y \rightarrow \mathbb R && \Cal B_{\text{nl}}(w, q; \boldsymbol{\mu}) = a_{\text{nl}}(z,q; \boldsymbol{\mu}) - c(v,q; \boldsymbol{\mu}) \hspace{1cm} &&& \forall w \in X, \; \forall q \in Y,\\
& F(\boldsymbol{\mu}) \in X\dual && \la F(\boldsymbol{\mu}), w \ra = m(y_d(\boldsymbol{\mu}),z; \boldsymbol{\mu}) \hspace{1cm} &&& \forall w \in X.
\end{align*}
Thanks to these definitions and taking into consideration the adjoint variable $p \in Y$ , we can define, as in the linear case, the functional \eqref{lin_func} and  the new Lagrangian functional $\Lg^{\text{nl}}: Y \times U \times Y \rightarrow \mathbb R$ as
\begin{equation}
\label{daje_nl}
\Lg^{\text{nl}} (y, u, p; \boldsymbol{\mu}) = \Cal J((y,u); \boldsymbol{\mu}) + \Cal B_{\text{nl}}((y,u), p; \boldsymbol{\mu}) - \la G(\boldsymbol{\mu}), p \ra.
\end{equation}
 In literature (see e.g. \cite{hinze2008optimization}) it is well known that solve the nonlinear minimisation problem \eqref{eq: controllo} is equivalent to solving the following system: 
 given $\boldsymbol{\mu} \in \Cal P$, find $(x(\boldsymbol{\mu}), p(\boldsymbol{\mu})) \in X \times 
Y$, where $x(\pmb \mu) = (y(\pmb \mu), u(\pmb \mu))$,  such that  
\begin{equation}
\label{ZETA_nl}
\begin{cases}
D_y\Lg^{\text{nl}}(y, u, p; \pmb \mu)[z] = 0 & \forall z \in Y,\\
D_u\Lg^{\text{nl}}(y, u, p; \pmb \mu)[v] = 0 & \forall v \in U,\\
D_p\Lg^{\text{nl}}(y, u, p; \pmb \mu) [q]= 0 & \forall q \in Y.\\
%\Cal A(x(\boldsymbol{\mu}), w; \boldsymbol{\mu}) + \Cal B_{\text{nl}}(w,p(\boldsymbol{\mu}); \boldsymbol{\mu}) = \la F(\boldsymbol{\mu}), w \ra & \forall w \in X, \\
%\Cal B_{\text{nl}} (x(\boldsymbol{\mu}), q; \boldsymbol{\mu}) = \la G(\boldsymbol{\mu}), q \ra & \forall q \in Y.
\end{cases}
\end{equation}

\no $D_y\Lg^{\text{nl}}(y, u, p; \pmb \mu), D_u\Lg^{\text{nl}}(y, u, p; \pmb \mu)$ and $D_p\Lg^{\text{nl}}(y, u, p; \pmb \mu)$ represent the differentiation of the Lagrangian functional \eqref{daje_nl} with respect to the state, the control ad the adjoint variable, respectively. 
\\Following the analogous argument proposed in section \ref{PROFE},  we consider $Y \disc$ and $U \disc$ as the Finite Element discretization for $Y$ and $U$, respectively. Also in this case, we define the discrete product space $X^{\Cal N} = Y \disc \times U \disc \subset X$ and the Galerkin Finite Element discretized version of the problem \eqref{ZETA_nl} as: given $\boldsymbol{\mu} \in \Cal P$, find $(x^{\Cal N}(\boldsymbol{\mu}), p^{\Cal N}(\boldsymbol{\mu})) \in X^{\Cal N} \times 
Y^ {\Cal N}$, where $x^\Cal N(\pmb \mu) = (y^\Cal N (\pmb \mu), u^\Cal N(\pmb \mu))$,  such that 

\begin{equation}
\label{FEOCP_nl}
\begin{cases}

D_{y\disc}\Lg^{\text{nl}}(y\disc, u\disc, p\disc; \pmb \mu) [z\disc]= 0 & \forall z\disc \in Y\disc,\\
D_{u\disc}\Lg^{\text{nl}}(y\disc, u\disc, p\disc; \pmb \mu) [v\disc]= 0 & \forall v\disc \in U\disc,\\
D_{p\disc}\Lg^{\text{nl}}(y\disc, u\disc, p\disc; \pmb \mu) [{q\disc}]= 0 & \forall q\disc \in Y\disc,\\

% \Cal A(x^{\Cal N}(\boldsymbol{\mu}),w\disc; \boldsymbol{\mu}) + \Cal B_{\text{nl}} (w\disc,p(\boldsymbol{\mu})\disc; \boldsymbol{\mu}) = \la F(\boldsymbol{\mu}), w\disc \ra & \forall v\disc \in X\disc, \\
%\Cal B_{\text{nl}}(x\disc(\boldsymbol{\mu}),q\disc; \boldsymbol{\mu})= \la G(\boldsymbol{\mu}), q\disc \ra & \forall q\disc \in Y\disc. 
\end{cases}
\end{equation}

\no where the discrete Lagrangian functional have been differentiated with respect to the discrete variables. Numerically, the discrete OCP($\pmb \mu$) problem (\ref{FEOCP_nl}) has been solved through Newton's methods. Since we are able to solve \emph{truth} problems, we can consider the construction of reduced basis functions exploiting the POD-Galerkin procedure illustrated in section \ref{PODsec}. Thanks to the POD-Galerkin algorithm we can reduce the Finite Element spaces $Y \disc$ and $U \disc$ and consider $Y_N$ and $U_N$, respectively. Then we can consider the space $X_N = Y_N \times U_N$ and the reduced nonlinear OCP($\pmb \mu$) to be solved is defined as:
given $\boldsymbol{\mu} \in \Cal P$, find $(x_N(\boldsymbol{\mu}), p_N(\boldsymbol{\mu})) \in X_N \times Y_N$, where $x_N(\pmb \mu) = (y_N (\pmb \mu), u_N(\pmb \mu))$, such that
\begin{equation}
\label{RBOCP_nl}
\begin{cases}

D_{y_N}\Lg^{\text{nl}}(y_N, u_N, p_N; \pmb \mu)[z_N] = 0 & \forall z_N \in Y_N,\\
D_{u_N}\Lg^{\text{nl}}(y_N, u_N, p_N; \pmb \mu)[v_N] = 0 & \forall v_N \in U_N,\\
D_{p_N}\Lg^{\text{nl}}(y_N, u_N, p_N; \pmb \mu) [q_N]= 0 & \forall q_N \in Y_N,\\
%\Cal A(x_N(\boldsymbol{\mu}),w_N; \boldsymbol{\mu}) + \Cal B_{\text{nl}}(w_N, %p_N(\boldsymbol{\mu}); \boldsymbol{\mu}) = \la F(\boldsymbol{\mu}), w_N \ra & \forall w_N \in X_N, \\
%\Cal B_{\text{nl}}(x_N(\boldsymbol{\mu}), q_N; \boldsymbol{\mu}) =  \la G(\boldsymbol{\mu}), q_N \ra & \forall q_N \in Y_N.
\end{cases}
\end{equation}

\no We recall that the affinity assumption must be verified in order to guarantee good performances of the ROM methods. 
Since the case of interest contains only quadratically nonlinear terms, the affinity assumption can be guaranteed by storing the appropriate
nonlinear terms in third order tensors.  In more general cases one can resort e.g. to empirical interpolation \cite{barrault2004empirical} and later variants.
We exploited the Newton's method to solve the reduced problem \eqref{RBOCP_nl}, as we did for the did for the discretized version \eqref{FEOCP_nl}.

\subsection{Reduced Basis Applied to a Nonlinear Ocean Circulation Solution Tracking}
\label{nl_sec}
We recall that the general Ocean circulation model governing large scale flow dynamics presented in \cite[Chapter 3]{cavallini2012quasi} is described by the nonlinear equation \eqref{QG}. The nonlinearity of the model is given by the expression \eqref{QGnl}. In this experiment, the physical parameter $\boldsymbol{\mu} = (\mu_1, \mu_2,\mu_3)$  takes values in the parameter space $\Cal P = [0.07^3,1]\times [10^{-4},1] \times [10^{-4}, 0.045^2]$. As in section \ref{OCEANO}, $\mu_1$ and $\mu_2$ represent the diffusive action of the Ocean, while $\mu_3$ is the parameter linked to its nonlinear dynamic. The range for the parameters ensures stability to the problem and allow us to treat a
\emph{moderate nonlinear} OCP($\pmb \mu$) governed by quasi-geostrophic equations: in this work we will only deal with this specific case and we will not treat \emph{highly nonlinear} problems, corresponding to lower values of $\mu_1, \mu_2$ and/or higher values of $\mu_3$. We underline that the steady quasi-geostrophic model is a very complex physical system: in order to have a complete description of the meaning of the model and of the role of all its components the reader is referred to section \ref{OCEANO} and to \cite[Chapter 3]{cavallini2012quasi}.
The solution tracking problem constrained to the nonlinear quasi-geostrophic equation is:

\begin{equation}
\label{pr_nl}
\begin{aligned}
& \underset{(\psi, u) \in V \times U}{\text{min }}
 J(\psi, u) =\underset{(\psi, u) \in V \times U}{\text{min }} \hspace{.2cm} \displaystyle \half \int_{\Omega} (\psi - \psi_d)^2 \; d \Omega + \displaystyle \alf \int_{\Omega} u^2 \; d \Omega \vspace{0.2cm}\\
& \qquad \text{such that } \hspace{0.5cm}
\begin{cases}
\displaystyle \frac{\partial \psi}{\partial x} = u
- \displaystyle \mu_1 \Delta \psi + \mu_2 \Delta ^2 \psi  - \mu_3 \Cal F(\psi, \Delta \psi)& \text{in } \Omega, \\
\psi = 0 & \text{on } \partial \Omega, \\
\Delta \psi = 0 & \text{on } \partial \Omega,
\end{cases}
\end{aligned}
\end{equation}
\no where $V$ and $U$ are suitable functional spaces, $\psi \in V$ is the state variable, $u \in U$ is the unknown wind action to be controlled and $\alpha$ is a penalization term. The aim of the OCP($\pmb \mu$) presented is the same of its linear version \eqref{nodec}: make our state solution the most similar to an already known desired state profile.
As in the linear case, we rewrite problem \eqref{pr_nl} as: 
\begin{equation}
\label{dec_nl}
\begin{aligned}
& \underset{((\psi ,q), u) \in Y \times U}{\text{min }}
 J((\psi, q), u) =\underset{((\psi ,q), u) \in Y \times U}{\text{min }}\hspace{.2cm} \displaystyle \half \int_{\Omega} (\psi - \psi_d)^2 \; d \Omega + \displaystyle \alf \int_{\Omega} u^2 \; d \Omega \vspace{0.2cm}\\
& \qquad \text{such that } \hspace{0.5cm}
\begin{cases}
q = \Delta \psi & \text{in } \Omega, \\
\displaystyle \frac{\partial \psi}{\partial x} = u 
- \displaystyle \mu_1 q + \mu_2 \Delta q -\mu_3 \Cal F(\psi, q)& \text{in } \Omega, \\
\psi = 0 & \text{on } \partial \Omega, \\
q = 0 & \text{on } \partial \Omega,
\end{cases}
\end{aligned}
\end{equation}
\no where the spaces are defined as $Y = H^1_0(\Omega)\times H^1_0(\Omega)$ and $U = L^2(\Omega)$. The weak formulation of the nonlinear state equation is represented as: 
\begin{equation}
\label{statodebole_nl}
a_{\text{nl}}((\psi, q), (\phi, r); \boldsymbol{\mu}) = c(u,(\phi, r)) \hspace{1cm} \forall \phi, r \in H^1_0(\Omega),
\end{equation}
\no where  $c: U \times Y \rightarrow \mathbb R$ is \eqref{stato2} and $a_{\text{nl}}: Y \times Y \rightarrow \mbb R$ is given by:
\begin{equation}
\displaystyle a_{\text{nl}}((\psi, q), (\phi, r); \boldsymbol{\mu})  = a((\psi, q), (\phi, r); \boldsymbol{\mu})
 - \mu_3 \int_{\Omega} \psi \Big ( \frac{\partial q}{\partial y}\frac{\partial r}{\partial x} -
\frac{\partial q}{\partial x}\frac{\partial r}{\partial y} \Big ) \; d\Omega,
\end{equation}

\no where $a((\psi, q), (\phi, r); \boldsymbol{\mu})$ is defined in \eqref{stato1}, and $G \in Y\dual$ is $G \equiv 0$. \vspace{.5cm} \\
In the following, we aim at recasting the problem \eqref{dec_nl} in the framework proposed in \eqref{ZETA_nl}. Let us define the product space 
$X = Y \times U$ and let $x = ((\psi, q),u)$ and $w = ((\chi, t),v)$ be two elements of $X$. Let us consider $s = (\phi, r) \in Y$ as our adjoint variable . Let us describe the following quantities:
$$
\begin{aligned}
& \Cal A : X \times X \rightarrow \mbb R \hspace{1cm} &&  \Cal A(x,w) = m((\psi, q), (\chi,t)) + \alpha n(u,v),  \\
& \Cal B_{\text{nl}} : X \times Y \rightarrow \mbb R \hspace{1cm} && \Cal B_{\text{nl}}(w, s; \boldsymbol{\mu})  = a_{\text{nl}}((\chi,t),(\phi, r), \boldsymbol{\mu}) - c(v, (\phi,r)), \\
& F: X \rightarrow \mathbb R \hspace{1cm} && \la F, w \ra = \int_{\Omega} \psi_d \chi \; d \Omega.\\
\end{aligned}
$$
As in the linear version, the bilinear forms $m : Y \times Y \rightarrow \mbb R$ and 
$n: U \times U \rightarrow \mbb R$ are 
$$
\displaystyle m((\psi, q), (\chi,t)) = \int_{\Omega} \psi \chi \; d \Omega \;\; \; \; \;\text{and}\; \;\;\;
\displaystyle n(u,v) = \int_{\Omega}uv \; d \Omega.
$$
\no
In order to solve the OCP($\pmb \mu$) governed by \eqref{dec_nl}, we define the Lagrangian functional \eqref{daje_nl} and we recall that $x = ((\psi, q),u), w = ((\chi, t),v)$. We also define a test function in the adjoint space $Y$ as $\eta = (\xi, \sigma)$. Then we solve the system \eqref{ZETA_nl} where: \vspace{-.3cm}
\begin{align*}
& D_{(\psi, q)}\Lg^{\text{nl}}((\psi, q), u, (\phi,r); \pmb \mu)[(\chi,t)] = m((\psi, q), (\chi,t)) 
- \int_{\Omega} \frac {\partial r }{\partial x}\chi \ \; d \Omega
\displaystyle + \mu_2  \int_{\Omega} \nabla t \cdot \nabla r \; d \Omega \; + \\
& \hspace{3cm}\qquad \qquad \displaystyle +\mu_1 \int_{\Omega} t r \; d \Omega 
 \displaystyle + \int_{\Omega} t \phi \; d \Omega
\displaystyle + \int_{\Omega} \nabla \chi \cdot \nabla \phi \; d \Omega \; + \vspace{0.2cm} \\
\tag{5.2.5}
& \hspace{3cm}\ \qquad \qquad \qquad \qquad \displaystyle - \mu_3
\int_{\Omega} \chi \Big ( \frac{\partial q}{\partial y}\frac{\partial r}{\partial x} -
\frac{\partial q}{\partial x}\frac{\partial r}{\partial y} \Big ) \; d\Omega \; +\\
& \hspace{3cm}\ \qquad \qquad \qquad \qquad \qquad \qquad \displaystyle - \mu_3
\int_{\Omega} t \Big ( \frac{\partial \psi}{\partial y}\frac{\partial r}{\partial y} -
\frac{\partial r}{\partial x}\frac{\partial \psi}{\partial y} \Big ) \; d\Omega -\int_{\Omega}u\chi - \la F, w \ra,\\
\displaystyle  &D_u \Lg^{\text{nl}}((\psi, q), u, (\phi,r); \pmb \mu)[v] = \alpha n(u,v) -   \int_{\Omega} v r \; d \Omega ,\\
\displaystyle & D_{(\phi, r)}\Lg^{\text{nl}}((\psi, q), u, (\phi, r); \pmb \mu)[ (\xi, \sigma)] =  \Cal B(x, \eta; \pmb \mu).
\end{align*}
We built the aggregated reduced spaces exploiting the partitioned POD-Galerkin algorithm for the five variables $\psi, q, u, \chi, t$ as proposed in section \ref{PODsec}, separately. As we did for the other numerical simulations, we underline that the affinity assumption is guaranteed. Indeed, with $Q_{\Cal A} = 1$, $Q_{\Cal B_{\text{nl}}} = 4$ and $Q_{F} =1$ the affine decomposition of the problem is given by

$$
\begin{aligned}
&\Theta^1_{\Cal A} = 1 \hspace{1cm} && \Cal A^1(x,w) = \Cal A(x,w), \\
&\Theta^1_{\Cal B_{\text{nl}}} = \mu_1 \hspace{1cm} && \Cal B^1_{\text{nl}}(x,s)= \displaystyle  \int_{\Omega} q r \; d \Omega,  \\
 &\Theta^2_{\Cal B_{\text{nl}}} = \mu_2 \hspace{1cm} && \Cal B^2_{\text{nl}}(x,s)  = \displaystyle \int_{\Omega} \nabla q \cdot \nabla r \; d \Omega,\\
 &\Theta^3_{\Cal B_{\text{nl}}} = 1 \hspace{1cm} && \Cal B^3_{\text{nl}}(x,s)  = \displaystyle\int_{\Omega} \frac {\partial \psi }{\partial x} r \; d \Omega
+ \displaystyle  \int_{\Omega} q \phi \; d \Omega
+ \int_{\Omega} \nabla \psi \cdot \nabla \phi \; d \Omega - \int_{\Omega}u r \; d \Omega,\\
&\Theta^4_{\Cal B_{\text{nl}}} = - \mu_3 \hspace{1cm} && \Cal B^4_{\text{nl}}(x,s)  = \int_{\Omega} \psi \Big ( \frac{\partial q}{\partial y}\frac{\partial r}{\partial x} -
\frac{\partial q}{\partial x}\frac{\partial r}{\partial y} \Big ) \; d\Omega,\\
&\Theta^1_F  = 1 \hspace{1cm} && \la F^1, w \ra = \la F, w \ra.\\
\end{aligned}
$$
%\no \textcolor{red}{\textbf{inizio aggiunta}}\\
%Since we are facing an optimal control problem governed by an elliptic coercive state equation, the well-posedness of the problem can be proved, in other words, conditions (i)-(iv) presented in section \ref{LQ} are fulfilled, as specified in \cite{negri2013reduced}. We focus our attention on the properties of the bilinear form $a(\cdot, \cdot)$. The properties needed by the forms $c(\cdot, \cdot)$, $n(\cdot, \cdot)$ and $n(\cdot, \cdot)$ to fulfill (i)-(iv) (they are presented in \cite{bochev2009least}) are trivial. \textcolor{blue}{\textbf{CREDETE CHE DOVREI ELENCARLE??}}
%Our purpose is to verify the hypotheses (i)-(iv) introduced in section \ref{conty}, in order to prove the well-posedness of the solution tracking optimal control problem. Let us %recall that inequality (D) holds for $X$ and $Q$. The following holds:
%\begin{enumerate}[(i)]
%\item the bilinear form  $\Cal A(\cdot, \cdot; \boldsymbol{\mu})$ is continuous over $X \times X$, since
%\begin{align*}
%|\Cal {A}(x, w)| & \leq \norm{\psi}_2\norm{\chi}_2 + \alpha \norm{u}_2\norm{v} \\
%			& \leq \max \{1,\alpha \} \norm{x}_X \norm{w}_X; \\
%\end{align*}
%\item
%\item the bilinear form  $\Cal B(\cdot, \cdot; \boldsymbol{\mu})$ is continuous over $X \times Q$, since
%\begin{align*}
%|\Cal {B}(x, s; \boldsymbol{\mu})| & \leq | a((\psi, q), (\phi, r); \boldsymbol{\mu})|+ |c(u, (\phi, r))| \\
%					& \leq \max \{ \mu_1, \mu_2, 1\}\norm{x}_X\norm{s}_Q;
%\end{align*}
%\item 
%\end{enumerate}
% \textcolor{red}{\textbf{fine aggiunta}}\\
\vspace{-.3cm}
\begin{table}[h]
\centering
\caption{\nl{Data of the numerical experiment: North Atlantic Ocean nonlinear solution tracking.}}
\label{antony_nl}
\begin{tabular}{ c | c }
\toprule
\textbf{Data} & \textbf{Values} \\
\midrule
 $(\mu_1, \mu_2, \mu_3,  \alpha)$ &  $(10^{-4}, 0.07^3, 0.045^2, 10^{-5})$ \\
\midrule
& FE solution of \\ $\psi_d$  & nonlinear quasi-geostrophic equation \\ & with $f = -\sin (\pi y)$\\ & and $\pmb \mu = (10^{-4}, 0.07^3, 0.07^2, 10^{-5})$ \\ \midrule
POD Training Set Dimension & 100 \\ 
\midrule
Basis Number $N$ & 25 \\
\midrule
Sampling Distribution & log-equispaced \\
\bottomrule
\end{tabular}
\end{table}

\no \\ Table \ref{antony_nl} shows all the features of this experiment. In Figure \ref{NLres} the desired solution profile $\psi_d$ is presented. Then the truth and the reduced solutions are plotted. It can be seen that the approximated solutions match. In the last plot of Figure \ref{NLres}, the pointwise error is shown: the maximum value reached is $9.6 \cdot 10^{-7}$ with \emph{basis number} $N = 25$.  We recall that $2 N = N_{\psi} = N_{q} = N_{\chi} = N_{t}$ and $N = N_u$. 
 We also analysed the ROM and the FE performances, in terms of: system dimension, cost functional optimal value, time of resolution (Table \ref{u3}). Table \ref{nlspeedup} presents the \textit{speed up} index with respect to the \textit{basis number} $N$. The results presented remark how the ROM approach could be a very suitable tool in order to solve quasi-geostrophic equations, most of all in their nonlinear version that describes more complicated, but more realistic, Ocean dynamics, like the movement of the flow stream towards North (Gulf Stream circulation). 
Figure \ref{err_nl} shows the error norm between FE and reduced variables as already indicated in footnote \ref{error} over a random testing set of $50$, obtaining a similar behavior with respect to the linear case.
\begin{table}[H]
\centering
\caption{\nl{ROM vs FE: North Atlantic Ocean nonlinear solution tracking.}}
\label{u3}
\begin{tabular}{c | c| cl}
\toprule
& FE & ROM \\
\midrule
System Dimension & $6490\times 6490$  & $450 \times 450$ \\
Optimal Cost Functional & $2.04347 \cdot 10^6$ & $2.04346 \cdot 10^6$\\
Time of Resolution & $7.21s$ & $5.11 \cdot 10^{-1}s$ \\
\bottomrule
\end{tabular}
\end{table}
\vspace{-.3cm}
\begin{figure}[h]

\hspace{-.8cm}\includegraphics[scale = 0.2]{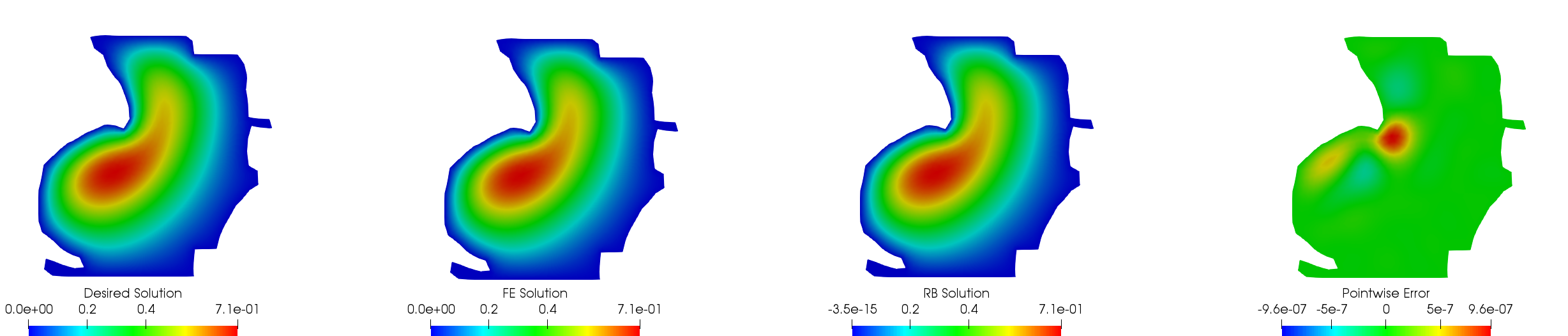}
%\hspace{-1.2cm}\includegraphics[scale = 0.25]{Atlantico/nl_FE_solution}
%\vspace{-.4cm}\hspace{-1.2cm}\includegraphics[scale = 0.25]{Atlantico/nl_RB_solution}\hspace{-1cm}\includegraphics[scale = 0.25]{Atlantico/nl_DIFF}
\caption{\nl{North Atlantic Ocean nonlinear solution tracking: results.}}
\label{NLres}
\end{figure}

\begin{figure}[H]
\centering
\hspace{-.5cm}
\includegraphics[scale = 0.28]{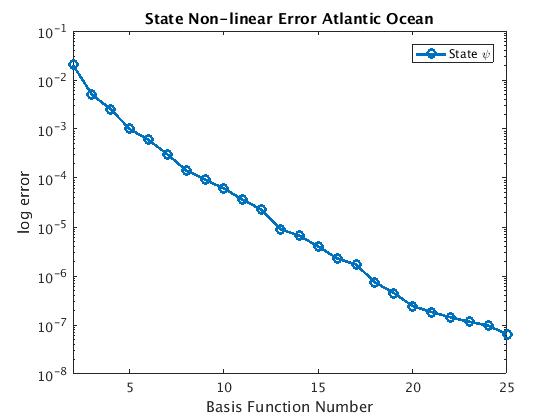}\includegraphics[scale = 0.28]{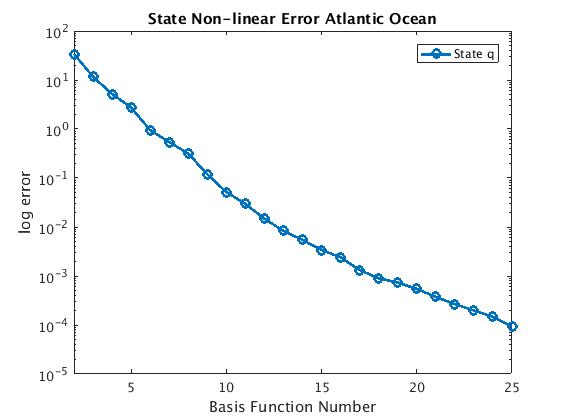}\includegraphics[scale = 0.28]{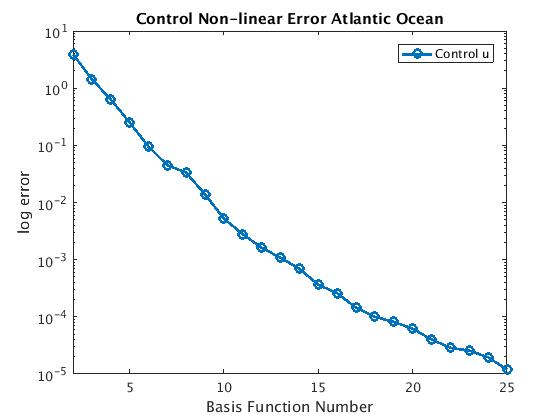}\\ 
\includegraphics[scale = 0.28]{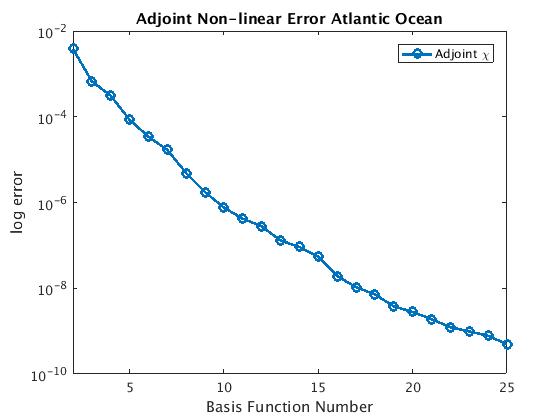}\includegraphics[scale = 0.28]{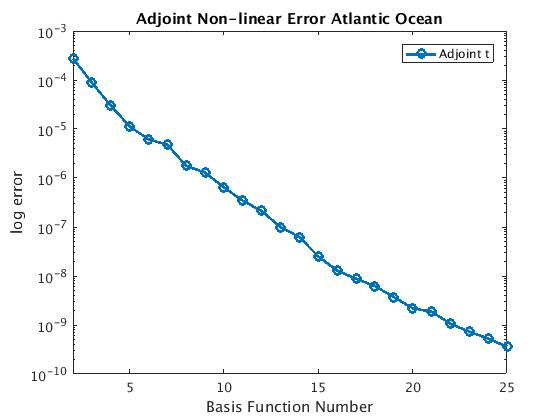}

\caption{\nl{North Atlantic Ocean nonlinear solution tracking: errors.}}
\label{err_nl}
\end{figure}

\begin{table}[H]
\centering
\caption{\nl{Speed up analysis: North Atlantic Ocean nonlinear solution tracking.}}
\label{nlspeedup}
\begin{tabular}{c | cccccccccccl}

\toprule

\nl{Basis Number } $\nl{N}$            &\nl{ 1} & \nl{5} & \nl{10} & \nl{15} &\nl{20}& \nl{25} \\
\nl{Speed up}     &       \nl{  10}& \nl{8} & \nl{7} &\nl{6} & \nl{5} & \nl{5} \\

\bottomrule
\end{tabular}

\end{table}

}

\section{Conclusions and Perspectives}
\label{conc}
In this work we have exploited reduced order methods in environmental parametrized optimal control problems dealing with marine sciences and engineering. We proposed two specific examples: one representing the ecological issue of pollutant control in a specific naturalist area, the Gulf of Trieste, Italy, the other one consisting in a large scale solution tracking OCP($\boldsymbol{\mu}$) governed by quasi-geostrophic equations. We showed how reduced order methods could be a  very useful tool in environmental sciences, like oceanography and ecology, where parametrized simulations are usually very demanding and costly. Reduced order methods are a suitable approach to face these issues. We have used a POD-Galerkin method for sampling and for the projection, by exploiting an aggregated space strategy. Reduced order methods performances have been compared to FE approximation, classically used to study these phenomena, in order to prove how convenient the reduced order approach could be in this particular field of applications: in the linear version of the North Atlantic problem, the reduced time of resolution decreases of one order of magnitude with respect to the \emph{full order} one, and the error of the state variable $\psi$ is negligible. To the best of our knowledge, the main novelty of this work is in the POD-Galerkin reduction of a solution tracking optimal control problem governed by quasi-geostrophic equations
\nl{in its linear and nonlinear version}. 
\\ 
\no Let us expose some improvements of this work, focusing on the optimal control problem governed by quasi-geostrophic equation. A possible development would involve a time dependent optimal control problem considering also the \nl{highly } nonlinear case. This kind of formulation is of the utmost importance in climatological applications, in order to forecast and predict possible scenarios in a reliable way. \nl{This complete model } will make the problems more and more realistic and suited to actual ecological and climatological challenges, as well as more and more computational demanding. In this sense, reduced order modelling appears, again, to be a suitable and versatile approach to be used. Time dependent nonlinear optimal control problems insert themselves in the framework of \textit{data assimilation techniques}, that, as briefly introduced in section \ref{OCEANO}, allow to modify the model in order to reach more reliable results in the forecasting applications, thanks to a solution tracking where the solution desired to be reached represents real experimental data.
\\ For all the examples presented, a further step could be the development of three-dimensional marine model that could take into consideration bathimetry effect. Finally, the problems could be inserted in a reduced order uncertainty quantification context (see e.g. \cite{Davide}), when it is not possible to assign specific values for the parameters by classical statistical methods.

\section*{Acknowledgements}
We acknowledge the support by European Union Funding for Research and Innovation -- Horizon 2020 Program -- in the framework of European Research Council Executive Agency: Consolidator Grant H2020 ERC CoG 2015 AROMA-CFD project 681447 ``Advanced Reduced Order Methods with Applications in Computational Fluid Dynamics''. We also acknowledge the INDAM-GNCS project ``Metodi numerici avanzati combinati con tecniche di riduzione computazionale per PDEs parametrizzate e applicazioni''.

%\nocite{*}

\bibliographystyle{plain} % "apa": author-year è lo stile e può essere modificato
%\bibliography{bibi}
\addcontentsline{toc}{chapter}{Bibliography}
\bibliography{BIB}
%\afterpage{\null\thispagestyle{empty}\clearpage}

\end{document}